\documentclass[12pt]{article}
\usepackage{amsmath,amsthm,amssymb,amsfonts}

\newtheorem{theorem}{Theorem}[section]
\newtheorem{proposition}[theorem]{Proposition}
\newtheorem{corollary}[theorem]{Corollary}
\newtheorem{lemma}[theorem]{Lemma}

\newtheorem{remark}[theorem]{Remark}

\begin{document}

\title{Torus $n$-Point Functions for $\mathbb{R}$-graded Vertex Operator
Superalgebras and Continuous Fermion Orbifolds }
\author{Geoffrey Mason\thanks{%
Partial support provided by NSF, NSA and the Committee on Research,
University of California, Santa Cruz} \\
Department of Mathematics, \\
University of California Santa Cruz, \\
CA 95064, U.S.A. \and Michael P. Tuite, Alexander Zuevsky\thanks{%
Supported by a Science Foundation Ireland Frontiers of Research Grant, and
by Max-Planck Institut f\"{u}r Mathematik, Bonn} \\
Department of Mathematical Physics, \\
National University of Ireland, \\
Galway, Ireland.\\
}
\maketitle

\begin{abstract}
We consider genus one $n$-point functions for a vertex operator superalgebra
with a real grading. We compute all $n$-point functions for rank one and
rank two fermion vertex operator superalgebras. In the rank two fermion
case, we obtain all orbifold $n$-point functions for a twisted module
associated with a continuous automorphism generated by a Heisenberg bosonic
state. The modular properties of these orbifold $n$-point functions are
given and we describe a generalization of Fay's trisecant identity for
elliptic functions.
\end{abstract}

\section{Introduction}

This paper is one of a series devoted to the study of $n$-point functions
for vertex operator algebras on Riemann surfaces of genus one, two and
higher \cite{T}, \cite{MT1}, \cite{MT2}, \cite{MT3}. One may define $n$%
-point functions at genus one following Zhu \cite{Z}, and use these
functions together with various sewing procedures to define $n$-point
functions at successively higher genera \cite{T}, \cite{MT2}, \cite{MT3}. In
this paper we consider the genus one $n$-point functions for a Vertex
Operator Superalgebra (VOSA) $V$ with a real grading (i.e. a chiral
fermionic conformal field theory). In particular, we compute all $n$-point
functions for rank one and rank two fermion VOSAs. In the latter case, we
consider $n$-point functions defined over\ an orbifold $g$-twisted module
for a continuous $V$ automorphism $g$ generated by a Heisenberg bosonic
state. We also consider the Heisenberg decomposition (or bosonization) of $V$
and recover elliptic versions of Fay's generalized trisecant identity
together with a new further generalization. The modular properties of the
continuous orbifold $n$-point functions are also described.

In his seminal paper, Zhu defined and developed a constructive theory of
torus $n$-point functions for a $\mathbb{Z}$-graded Vertex Operator Algebra
(VOA) and its modules \cite{Z}. In particular, he described various
recursion formulae where, for example, an $n$-point function is expanded in
terms of $n-1$-point functions and naturally occurring Weierstrass elliptic
(and quasi-elliptic) functions. Indeed, one can prove the analytic, elliptic
and modular properties of $n$-point functions for many VOAs from these
recursion formulae (op.cit.). This technique has since been extended to
include orbifold VOAs with a $g$-twisted module for a finite order
automorphism $g$ \cite{DLM1}, $\frac{1}{2}\mathbb{Z}$-graded VOSAs \cite{DZ1}
and $\mathbb{Z}$-graded VOSAs \cite{DZ2}. Here we consider a further
generalization to obtain recursion formulae for torus $n$-point functions
for an $\mathbb{R}$-graded VOSA. We consider $n$-point functions defined as
the supertrace over the product of various vertex operators together with a
general element of the automorphism group of the VOSA. The resulting
recursion formula is expressed in terms of natural "twisted" Weierstrass
elliptic functions periodic up to arbitrary multipliers in $U(1)$. Such
elliptic functions already appear in ref. \cite{DLM1} for multipliers of
finite order. Here, we give a detailed description of twisted Weierstrass
elliptic functions (and associated twisted Eisenstein series) for general $%
U(1)$ multipliers generalizing many results of the classical theory of
elliptic functions.

We consider two applications of the Zhu recursion formula. The first example
is that of the rank one $\frac{1}{2}\mathbb{Z}$-graded fermion VOSA. In this
case, all $n$-point functions can be computed in terms of a single
generating function. In particular, we obtain expressions for these $n$%
-point functions in a natural Fock basis in terms of the Pfaffian of an
appropriate block matrix. The second example is that of the rank two fermion
VOSA. As is well known, this VOSA contains a Heisenberg vector which
generates a continuous automorphism $g$ and for which a $g$-twisted module
can be constructed \cite{Li}. The Heisenberg vector can also be employed to
define a "shifted" Virasoro with real grading \cite{MN2}, \cite{DM}. We
demonstrate a general relationship between the $n$-point functions for
orbifold $g$-twisted modules and the shifted VOSA. We next apply the
recursion formula for $\mathbb{R}$-graded VOSAs in order to obtain all
continuous orbifold $n$-point functions. These are expressible in terms of
determinants of appropriate block matrices in a natural Fock basis and can
again be obtained from a single generating function. Decomposing the rank
two fermion VOSA into Heisenberg irreducible modules as a bosonic $\mathbb{Z}
$-lattice VOSA (i.e. bosonization) we may employ results of ref. \cite{MT2}
to find alternative expressions for the $n$-point functions. In particular
the generating function is expressible in terms of theta functions and the
genus one prime form and we thus recover Fay's generalized trisecant
identity for elliptic functions. We also prove a further generalization of
Fay's trisecant identity based on the $n$-point function for $n$ lattice
vectors. The paper concludes with a determination of the modular
transformation properties for all rank two continuous orbifold $n$-point
functions generalizing Zhu's results for $C_{2}$-cofinite VOAs \cite{Z}.

The study of $n$-point functions has a long history in the theoretical
physics literature and we recover a number of well known physics results
here. Thus the Pfaffian and determinant formulas for the rank one and two
fermion generating functions and the relationship between Fay's generalized
trisecant identity have previously appeared in physics \cite{R1, R2, EO, RS,
FMS, P}. However, it is important to emphasize that our approach is
constructively based on the properties of a VOSA and that a rigorous and
complete description of these $n$-point functions has been lacking until
now. Thus, for example, no assumption is made about the local analytic
properties of $n$-point functions as would normally be the case in physics.
Similarly, other pure mathematical algebraic geometric approaches to $n$%
-point functions are based on an assumed local analytic structure \cite{TUY}%
. Finally, apart from the intrinsic benefits of this rigorous approach, it
is important to obtain a complete description of these $n$-point functions
as the building blocks used in the construction of higher genus partition
and $n$-point functions \cite{T}, \cite{MT2}, \cite{MT3}.

The paper is organized as follows. We begin in Section 2 with a review of
classical Weierstrass elliptic functions and Eisenstein series. We introduce
twisted Weierstrass functions which are periodic up to arbitrary elements of 
$U(1)$. We describe various expansions of these twisted functions, introduce
twisted Eisenstein series and determine their modular properties. 
Section 3 contains one of the central results of this paper. We begin with
the defining properties of an $\mathbb{R}$-graded VOSA $V$. We define $n$%
-point functions as a supertrace over a $V$-module and describe some general
properties. 
We then formulate a generalization of Zhu's recursion formula \cite{Z} to an 
$\mathbb{R}$-graded VOSA module making use of the twisted Weierstrass and
Eisenstein series. 
Section 4 contains a discussion of a VOSA containing a Heisenberg vector. We
prove the general relationship between the $n$-point functions for a VOSA
with a Heisenberg shifted Virasoro vector and $g$-twisted $n$-point
functions where $g$ is generated by the Heisenberg vector. 
In section 5 we apply the results of Section 3 to a rank one fermion VOSA.
In particular, we compute all $n$-point functions in terms of a generating
function given by a particular $n$-point function. We also discuss $n$-point
functions for a fermion number-twisted module. 
Section 6 contains a description of a rank two fermion VOSA. 
We make use of the results of Section 3 and Section 4 to compute all $n$%
-point functions for a $g$-twisted module where $g$ is generated by a
Heisenberg vector by means of a generating function. We next discuss the
Heisenberg decomposition of this rank two theory - the bosonized theory. In
particular, we derive an expression for the rank two generating function in
terms of $\theta $-functions and prime forms related to Fay's generalized
trisecant identity for elliptic functions. A further generalization for
Fay's trisecant identity for elliptic functions is also discussed. 
Finally, we discuss the modular properties of all $n$-point functions for
the rank two fermion VOSA. Properties of supertraces are recalled in the
Appendix.


We collect here notation for some of the more frequently occurring functions
and symbols that will play a role in our work. $\mathbb{Z}$ is the set of
integers, $\mathbb{R}$ the real numbers, $\mathbb{C}$ the complex numbers, $%
\mathbb{H}$ the complex upper-half plane. We will always take $\tau $ to lie
in $\mathbb{H}$, and $z$ will lie in $\mathbb{C}$ unless otherwise noted.
For a symbol $z\,$\ we set $q_{z}=\exp (z)$, in particular, $q=q_{2\pi i\tau
}=\exp (2\pi i\tau )$.

\section{\label{Section_Elliptic}Some Elliptic Function Theory}

\subsection{Classical Elliptic Functions}

We discuss a number of modular and elliptic-type functions that we will
need. We begin with some standard elliptic functions \cite{La}. The
Weierstrass $\wp$-function periodic in $z$ with periods $2\pi i$ and $2\pi
i\tau $ is 
\begin{eqnarray}
\wp (z,\tau ) &=&\frac{1}{z^{2}}+\sum_{\substack{ m,n\in \mathbb{Z}  \\ %
(m,n)\neq (0,0)}}[\frac{1}{(z-\omega _{m,n})^{2}}-\frac{1}{\omega _{m,n}^{2}}%
].  \label{Weier} \\
&=&\frac{1}{z^{2}}+\sum\limits_{n\geq 4,n\text{ even}}(n-1)E_{n}(\tau
)z^{n-2},  \label{Weier_En}
\end{eqnarray}%
for $(z,\tau )\in \mathbb{C}\times \mathbb{H}$ with $\omega _{m,n}=2\pi
i(m\tau +n)$. Here, $E_{n}(\tau )$ is equal to $0$ for $n$ odd, and for $n$
even is the Eisenstein series \cite{Se} 
\begin{equation}
E_{n}(\tau )=-\frac{B_{n}(0)}{n!}+\frac{2}{(n-1)!}\sum\limits_{r\geq 1}\frac{%
r^{n-1}q^{r}}{1-q^{r}},  \label{Eisen}
\end{equation}%
where $B_{n}(0)$ is the $n$th Bernoulli number (see (\ref{Bernoulli})
below). If $n\geq 4$ then $E_{n}(\tau )$ is a holomorphic modular form of
weight $n$ on $SL(2,\mathbb{Z})$. That is, it satisfies 
\begin{equation}
E_{n}(\gamma .\tau )=(c\tau +d)^{n}E_{n}(\tau ),  \label{gammaEn}
\end{equation}%
for all $\gamma =\left( 
\begin{array}{cc}
a & b \\ 
c & d%
\end{array}%
\right) \in SL(2,\mathbb{Z})$, where we use the standard notation

\begin{equation}
\gamma .\tau =\frac{a\tau +b}{c\tau +d}.  \label{mobiusaction}
\end{equation}%
On the other hand, $E_{2}(\tau )$ is a quasimodular form \cite{KZ} having
the exceptional transformation law 
\begin{equation}
E_{2}(\gamma .\tau )=(c\tau +d)^{2}E_{2}(\tau )-\frac{c(c\tau +d)}{2\pi i}.
\label{gammaE2}
\end{equation}

We define $P_{k}(z,\tau )$ for $k\geq 1$ by 
\begin{equation}
P_{k}(z,\tau )=\frac{(-1)^{k-1}}{(k-1)!}\frac{d^{k-1}}{dz^{k-1}}P_{1}(z,\tau
)=\frac{1}{z^{k}}+(-1)^{k}\sum\limits_{n\geq k}\binom{n-1}{k-1}E_{n}(\tau
)z^{n-k}.  \label{Pkdef}
\end{equation}%
Then $P_{2}(z,\tau )=\wp (z,\tau )+E_{2}(\tau )$ whereas $P_{1}-zE_{2}$ is
the classical Weierstrass zeta function. $P_{k}$ has periodicities 
\begin{eqnarray}
P_{k}(z+2\pi i,\tau ) &=&P_{k}(z,\tau ),  \notag \\
P_{k}(z+2\pi i\tau ,\tau ) &=&P_{k}(z,\tau )-\delta _{k1}.  \label{Pkperiod}
\end{eqnarray}

We define the elliptic prime form $K(z,\tau )\,$\ by \cite{Mu} 
\begin{equation}
K(z,\tau )=\exp (-P_{0}(z,\tau )),  \label{Primeform}
\end{equation}%
where 
\begin{equation}
P_{0}(z,\tau )=-\log (z)+\sum_{k\geq 2}\frac{1}{k}E_{k}(\tau )z^{k},
\label{P0}
\end{equation}%
so that 
\begin{equation}
P_{1}(z,\tau )=-\frac{d}{dz}P_{0}(z,\tau )=\frac{1}{z}-\sum_{k\geq
2}E_{k}(\tau )z^{k-1}.  \label{P1_dlogKdz}
\end{equation}%
$K(z,\tau )$ has periodicities 
\begin{eqnarray}
K(z+2\pi i,\tau ) &=&-K(z,\tau ),  \notag \\
K(z+2\pi i\tau ,\tau ) &=&-q_{z}^{-1}q^{-1/2}K(z,\tau ).
\label{Kquasiperiod}
\end{eqnarray}

We define the standard Jacobi theta function by\footnote{%
Note that the $z$ dependence of the theta function is chosen so that the
periods are $2\pi i$ and $2\pi i\tau $ rather than the standard periods of $%
1 $ and $\tau $.} e.g. \cite{FK} 
\begin{equation}
\vartheta \left[ 
\begin{array}{c}
a \\ 
b%
\end{array}%
\right] (z,\tau )=\sum\limits_{n\in \mathbb{Z}}\exp [i\pi (n+a)^{2}\tau
+(n+a)(z+2\pi ib)],  \label{thetaab}
\end{equation}%
with periodicities%
\begin{eqnarray}
\vartheta \left[ 
\begin{array}{c}
a \\ 
b%
\end{array}%
\right] (z+2\pi i,\tau ) &=&e^{2\pi ia}\vartheta \left[ 
\begin{array}{c}
a \\ 
b%
\end{array}%
\right] (z,\tau ),  \label{Thetaperiod_2pi_i} \\
\vartheta \left[ 
\begin{array}{c}
a \\ 
b%
\end{array}%
\right] (z+2\pi i\tau ,\tau ) &=&e^{-2\pi ib}q_{z}^{-1}q^{-1/2}\vartheta %
\left[ 
\begin{array}{c}
a \\ 
b%
\end{array}%
\right] (z,\tau ).  \label{Thetaperiod_2Pi_i_tau}
\end{eqnarray}%
We also note the modular transformation properties under the action of the
standard generators $S=\left( 
\begin{array}{cc}
0 & 1 \\ 
-1 & 0%
\end{array}%
\right) $ and $T=\left( 
\begin{array}{cc}
1 & 1 \\ 
0 & 1%
\end{array}%
\right) $ of $SL(2,\mathbb{Z})$ (with relations $(ST)^{3}=-S^{2}=1$) 
\begin{eqnarray}
\vartheta \left[ 
\begin{array}{c}
a \\ 
b%
\end{array}%
\right] (z,\tau +1) &=&e^{-i\pi a(a+1)}\vartheta \left[ 
\begin{array}{c}
a \\ 
b+a+\frac{1}{2}%
\end{array}%
\right] (z,\tau ),  \label{S_Theta} \\
\vartheta \left[ 
\begin{array}{c}
a \\ 
b%
\end{array}%
\right] (-\frac{z}{\tau },-\frac{1}{\tau }) &=&(-i\tau )^{1/2}e^{2\pi
iab}e^{-iz^{2}/4\pi \tau }\vartheta \left[ 
\begin{array}{c}
-b \\ 
a%
\end{array}%
\right] (z,\tau ).  \label{T_Theta}
\end{eqnarray}

$K(z,\tau )$ can be expressed in terms of half integral theta functions as 
\begin{equation}
K(z,\tau )=\frac{\vartheta \left[ 
\begin{array}{c}
\frac{1}{2} \\ 
\frac{1}{2}%
\end{array}%
\right] (z,\tau )}{\frac{d}{dz}\vartheta \left[ 
\begin{array}{c}
\frac{1}{2} \\ 
\frac{1}{2}%
\end{array}%
\right] (0,\tau )}=\frac{-i}{\eta (\tau )^{3}}\vartheta \left[ 
\begin{array}{c}
\frac{1}{2} \\ 
\frac{1}{2}%
\end{array}%
\right] (z,\tau ).  \label{Ktheta}
\end{equation}
where the Dedekind eta-function is defined by 
\begin{equation}
\eta (\tau )=q^{1/24}\prod_{n=1}^{\infty }(1-q^{n}).  \label{etafun}
\end{equation}

\subsection{\label{Section_Twisted}Twisted Elliptic Functions}

Let $(\theta ,\phi )\in U(1)\times U(1)$ denote a pair of modulus one
complex parameters with $\phi =\exp (2\pi i\lambda )$ for $0\leq \lambda <1$%
. For $z\in \mathbb{C}$ and $\tau \in \mathbb{H}$ we define "twisted"
Weierstrass functions for $k\geq 1$ as follows:%
\begin{equation}
P_{k}\left[ 
\begin{array}{c}
\theta \\ 
\phi%
\end{array}%
\right] (z,\tau )=\frac{(-1)^{k}}{(k-1)!}\sum\limits_{n\in \mathbb{Z}%
+\lambda }^{\prime }\frac{n^{k-1}q_{z}^{n}}{1-\theta ^{-1}q^{n}},
\label{Pkuv}
\end{equation}%
for $q=q_{2\pi i\tau }$ where $\sum\limits^{\prime }$ means we omit $n=0$ if 
$(\theta ,\phi )=(1,1)$.

\begin{remark}
\label{Remark_Pk}

\ (i) (\ref{Pkuv}) was introduced in \cite{DLM1} for rational $\lambda$,
where it was denoted by $P_{k}(\phi ,\theta ^{-1},z,\tau )$. The alternative
definition and notation used here is motivated by the modular and
periodicity properties shown below and by the column vector notation for
theta series.

(ii) (\ref{Pkuv}) converges absolutely and uniformly on compact subsets of
the domain $\left\vert q\right\vert <\left\vert q_{z}\right\vert <1$ \cite%
{DLM1}.

(iii) For $k\geq 1$,%
\begin{equation}
P_{k}\left[ {%
\begin{array}{c}
\theta \\ 
\phi%
\end{array}%
}\right] (z,\tau )=\frac{(-1)^{k-1}}{(k-1)!}\frac{d^{k-1}}{dz^{k-1}}P_{1}%
\left[ {%
\begin{array}{c}
\theta \\ 
\phi%
\end{array}%
}\right] (z,\tau ).  \label{Pkplus1}
\end{equation}
\end{remark}

We now develop twisted versions of the standard results for the classical
Weierstrass $\wp$-function reviewed above. A number of similar results
appear in \cite{DLM1}. However, the cases $k=1,2$ are treated separately
there and only for rational $\lambda $ i.e. $\phi ^{N}=1$ for some positive
integer $N$. The most canonical derivation of the periodic and modular
properties of (\ref{Pkuv}) for general $\lambda $ follow from the following
theorem:

\begin{theorem}
\label{Theorem_Pk_doublesum} For $\left\vert q\right\vert <\left\vert
q_{z}\right\vert <1$ and $\phi \neq 1$, 
\begin{equation}
P_{k}\left[ {%
\begin{array}{c}
\theta \\ 
\phi%
\end{array}%
}\right] (z,\tau )=\sum_{m\in \mathbb{Z}}\theta ^{m}\left[ \sum_{n\in 
\mathbb{Z}}\frac{\phi ^{n}}{\left( z-\omega _{m,n}\right) ^{k}}\right] ,
\label{P1_double1}
\end{equation}%
whereas for $\theta \neq 1$, 
\begin{equation}
P_{k}\left[ {%
\begin{array}{c}
\theta \\ 
\phi%
\end{array}%
}\right] (z,\tau )=\sum_{n\in \mathbb{Z}}\phi ^{n}\left[ \sum_{m\in \mathbb{Z%
}}\frac{\theta ^{m}}{\left( z-\omega _{m,n}\right) ^{k}}\right] .
\label{P1_double2}
\end{equation}
\end{theorem}

\begin{remark}
\label{Remark_Pk_doublesum} When \emph{both} $\theta \neq 1$ and $\phi \neq
1 $ then the double sums (\ref{P1_double1}) and (\ref{P1_double2}) are
equal. For $k\geq 3$, they are absolutely convergent and equal for all $%
(\theta ,\phi )$.
\end{remark}

In order to prove Theorem \ref{Theorem_Pk_doublesum} it is useful to define
the following convergent sum 
\begin{equation}
S(x,\phi )=\sum_{n\in \mathbb{Z}}\frac{\phi ^{n}}{x-2\pi in}.
\label{S(x,phi)}
\end{equation}%
Clearly%
\begin{eqnarray}
S(x+2\pi i,\phi ) &=&\phi S(x,\phi ),  \label{Speriod} \\
S(x,\phi ) &=&-S(-x,\phi ^{-1}).  \label{Sreflect}
\end{eqnarray}%
We then have:

\begin{lemma}
\label{Lemma_Sum_S(x,phi)} \ For $\phi =\exp (2\pi i\lambda )$ with $0\leq
\lambda <1$ we have 
\begin{equation}
S(x,\phi )=\frac{1}{2}\delta _{\lambda ,0}+\frac{q_{x}^{\lambda }}{q_{x}-1}.
\label{Sum_S(x,phi)}
\end{equation}
\end{lemma}

\noindent \textbf{Proof.} Both $S(x,\phi )$ and $q_{x}^{\lambda
}(q_{x}-1)^{-1}$ have simple poles at $x=2\pi in$ with residue $\phi ^{n}$
for all $n\in \mathbb{Z} $. Furthermore, $q_{x}^{\lambda }(q_{x}-1)^{-1}$ is
regular at the point at infinity for $0\leq \lambda <1$. Thus $S(x,\phi
)-q_{x}^{\lambda }(q_{x}-1)^{-1}$ is constant which from (\ref{Speriod}) and
(\ref{Sreflect}) must be given by $\frac{1}{2}\delta _{\lambda ,0}$. $%
\square $

\medskip We first prove Theorem \ref{Theorem_Pk_doublesum} for the case $k=1$
and $\phi \neq 1$ (i.e. $0<\lambda <1$). The double sum (\ref{P1_double1})
is 
\begin{equation*}
\sum_{m\in \mathbb{Z}}\theta ^{m}\left[ \sum_{n\in \mathbb{Z}}\frac{\phi ^{n}%
}{z-\omega _{m,n}}\right] =\sum_{m\in \mathbb{Z}}\theta ^{m}S(x_{m},\phi
)=\sum_{m\in \mathbb{Z}}\theta ^{m}\frac{q_{x_{m}}^{\lambda }}{q_{x_{m}}-1},
\end{equation*}%
using Lemma \ref{Lemma_Sum_S(x,phi)} for $x_{m}=z-2\pi im\tau $ with $%
q_{x_{m}}=q_{z}q^{-m}$. Since $\left\vert q\right\vert <\left\vert
q_{z}\right\vert <1$ we find for $m>0$ that $\left\vert q_{x_{m}}\right\vert
>1$ and hence 
\begin{equation*}
\frac{q_{x_{m}}^{\lambda }}{q_{x_{m}}-1}=\sum\limits_{r\leq
-1}q_{z}^{r+\lambda }(q^{-r-\lambda })^{m}.
\end{equation*}%
Since $\left\vert \theta q^{-r-\lambda }\right\vert <1$ for $r\leq -1$ we
obtain 
\begin{eqnarray*}
\sum_{m>0}\theta ^{m}\left[ \sum_{n\in \mathbb{Z}}\frac{\phi ^{n}}{z-\omega
_{m,n}}\right] &=&\sum\limits_{r\leq -1}q_{z}^{r+\lambda }\sum_{m>0}(\theta
q^{-r-\lambda })^{m} \\
&=&-\sum\limits_{r\leq -1}\frac{q_{z}^{r+\lambda }}{1-\theta
^{-1}q^{r+\lambda }}.
\end{eqnarray*}%
Similarly for $m\leq 0$ we have $\left\vert q_{x_{m}}\right\vert <1$, so
that 
\begin{equation*}
\frac{q_{x_{m}}^{\lambda }}{q_{x_{m}}-1}=-\sum\limits_{r\geq
0}q_{z}^{r+\lambda }(q^{-r-\lambda })^{m}.
\end{equation*}%
Hence since $\left\vert \theta q^{r+\lambda }\right\vert <1$ for $r\geq 0$
we find 
\begin{equation*}
\sum_{m\leq 0}\theta ^{m}\left[ \sum_{n\in \mathbb{Z}}\frac{\phi ^{n}}{%
z-\omega _{m,n}}\right] =-\sum\limits_{r\geq 0}\frac{q_{z}^{r+\lambda }}{%
1-\theta ^{-1}q^{r+\lambda }}.
\end{equation*}%
Altogether we obtain%
\begin{equation*}
\sum_{m\in \mathbb{Z}}\theta ^{m}\left[ \sum_{n\in \mathbb{Z}}\frac{\phi ^{n}%
}{z-\omega _{m,n}}\right] =-\sum\limits_{r\in \mathbb{Z}}\frac{%
q_{z}^{r+\lambda }}{1-\theta ^{-1}q^{r+\lambda }}=P_{1}\left[ {%
\begin{array}{c}
\theta \\ 
\phi%
\end{array}%
}\right] (z,\tau ),
\end{equation*}%
proving (\ref{P1_double1}) for $k=1$. The result for $k\geq 2$ follows after
applying (\ref{Pkplus1}).

In order to prove (\ref{P1_double2}) it is useful to first consider the
following double sum for $\phi \neq 1$ 
\begin{equation*}
A\left[ {%
\begin{array}{c}
\theta \\ 
\phi%
\end{array}%
}\right] (z,\tau )=\sum_{m\in \mathbb{Z}}\theta ^{m}\left[ \sum_{n\in 
\mathbb{Z}}\phi ^{n}\left( \frac{1}{z-\omega _{m,n}}-\frac{2}{z-\omega
_{m,n-1}}+\frac{1}{z-\omega _{m,n-2}}\right) \right] .
\end{equation*}%
By (\ref{Speriod}) we find%
\begin{eqnarray}
A\left[ {%
\begin{array}{c}
\theta \\ 
\phi%
\end{array}%
}\right] (z,\tau ) &=&\sum_{m\in \mathbb{Z}}\theta ^{m}\left[ S(x_{m},\phi
)-2S(x_{m}+2\pi i,\phi )+S(x_{m}+4\pi i,\phi )\right]  \notag \\
&=&(1-\phi )^{2}P_{1}\left[ {%
\begin{array}{c}
\theta \\ 
\phi%
\end{array}%
}\right] (z,\tau ).  \label{Asum1}
\end{eqnarray}%
On the other hand, we have%
\begin{equation*}
A\left[ {%
\begin{array}{c}
\theta \\ 
\phi%
\end{array}%
}\right] (z,\tau )=\sum_{m\in \mathbb{Z}}\theta ^{m}\left[ \sum_{n\in 
\mathbb{Z}}\phi ^{n}\frac{-8\pi ^{2}}{(z-\omega _{m,n})(z-\omega
_{m,n-1})(z-\omega _{m,n-2})}\right] .
\end{equation*}%
This sum is absolutely convergent since the summand is $O(\left\vert \omega
_{m,n}\right\vert ^{-3})$ for $\left\vert m\right\vert $, $\left\vert
n\right\vert $ large. We may thus interchange the order of summation to find
that, on relabelling, $A\left[ {%
\begin{array}{c}
\theta \\ 
\phi%
\end{array}%
}\right] (z,\tau )$ becomes 
\begin{eqnarray}
&&\sum_{m\in \mathbb{Z}}\phi ^{m}\left[ \sum_{n\in \mathbb{Z}}\theta
^{-n}\left( \frac{1}{z-\omega _{-n,m}}-\frac{2}{z-\omega _{-n,m-1}}+\frac{1}{%
z-\omega _{-n,m-2}}\right) \right]  \notag \\
&=&(-\frac{1}{\tau }) \sum_{m\in \mathbb{Z}}\phi ^{m}\left[ \sum_{n\in 
\mathbb{Z}}\theta ^{-n}\left( \frac{1}{z^{\prime }-\omega _{m,n}^{\prime }}-%
\frac{2}{z^{\prime }-\omega _{m-1,n}^{\prime }}+\frac{1}{z^{\prime }-\omega
_{m-2,n}^{\prime }}\right) \right]  \notag \\
&=&(-\frac{1}{\tau }) \sum_{m\in \mathbb{Z}}\phi ^{m}\left[ S(x_{m}^{\prime
},\theta ^{-1})-2S(x_{m-1}^{\prime },\theta ^{-1})+S(x_{m-2}^{\prime
},\theta ^{-1})\right],  \label{Asum2}
\end{eqnarray}%
where 
\begin{equation}
z^{\prime }=-\frac{z}{\tau },\quad \tau ^{\prime }=-\frac{1}{\tau },\quad
\omega _{m,n}^{\prime }=2\pi i(m\tau ^{\prime }+n),\quad x_{m}^{\prime
}=z^{\prime }-2\pi im\tau ^{\prime }.  \label{ztauprime}
\end{equation}%
Applying Lemma \ref{Lemma_Sum_S(x,phi)} with $\theta =\exp (-2\pi i\mu )$
for $0\leq \mu <1$, it follows that 
\begin{eqnarray*}
&&S(x_{m}^{\prime },\theta ^{-1})-2S(x_{m-1}^{\prime },\theta
^{-1})+S(x_{m-2}^{\prime },\theta ^{-1}) \\
&=&(1-2+1).\frac{1}{2}\delta _{\mu ,0}+\frac{q_{x_{m}^{\prime }}^{\mu }}{%
q_{x_{m}^{\prime }}-1}-2\frac{q_{x_{m-1}^{\prime }}^{\mu }}{%
q_{x_{m-1}^{\prime }}-1}+\frac{q_{x_{m-2}^{\prime }}^{\mu }}{%
q_{x_{m-2}^{\prime }}-1}.
\end{eqnarray*}%
We may next repeat the arguments above leading to (\ref{Asum1}) to find that
(\ref{Asum2}) becomes 
\begin{equation*}
A\left[ {%
\begin{array}{c}
\theta \\ 
\phi%
\end{array}%
}\right] (z,\tau )=(-\frac{1}{\tau })(1-\phi )^{2}P_{1}\left[ {%
\begin{array}{c}
\phi \\ 
\theta ^{-1}%
\end{array}%
}\right] (-\frac{z}{\tau },-\frac{1}{\tau }).
\end{equation*}%
Comparing to (\ref{Asum1}), we find that for $\phi \neq 1$ 
\begin{equation}
P_{1}\left[ {%
\begin{array}{c}
\theta \\ 
\phi%
\end{array}%
}\right] (z,\tau )=(-\frac{1}{\tau })P_{1}\left[ {%
\begin{array}{c}
\phi \\ 
\theta ^{-1}%
\end{array}%
}\right] (-\frac{z}{\tau },-\frac{1}{\tau }).  \label{P1mod_S}
\end{equation}%
Considering this identity for $(z^{\prime },\tau ^{\prime })$ of (\ref%
{ztauprime}) and using 
\begin{equation}
P_{1}\left[ {%
\begin{array}{c}
\theta \\ 
\phi%
\end{array}%
}\right] (z,\tau )=-P_{1}\left[ {%
\begin{array}{c}
\theta ^{-1} \\ 
\phi ^{-1}%
\end{array}%
}\right] (-z,\tau ),  \label{P1minusz}
\end{equation}%
(which follows from (\ref{P1_double1})) it is clear that (\ref{P1mod_S})
holds for all $(\theta ,\phi )\neq (1,1)$.

We may use (\ref{P1mod_S}) to prove (\ref{P1_double2}) of Theorem \ref%
{Theorem_Pk_doublesum} in the case $k=1$. The double sum of (\ref{P1_double2}%
) becomes on relabelling 
\begin{eqnarray*}
\sum_{m\in \mathbb{Z}}\phi ^{m}\left[ \sum_{n\in \mathbb{Z}}\frac{\theta
^{-n}}{z-\omega _{-n,m}}\right] &=&(-\frac{1}{\tau })P_{1}\left[ {%
\begin{array}{c}
\phi \\ 
\theta ^{-1}%
\end{array}%
}\right] (-\frac{z}{\tau },-\frac{1}{\tau }) \\
&=&P_{1}\left[ {%
\begin{array}{c}
\theta \\ 
\phi%
\end{array}%
}\right] (z,\tau ).
\end{eqnarray*}%
The general result for $k\geq 2$ follows from (\ref{Pkplus1}). $\square $

\medskip Periodicity and modular properties now follow from Theorem \ref%
{Theorem_Pk_doublesum}. Thus we have

\begin{lemma}
\label{Lemma_periodicity} \ For $(\theta ,\phi )\neq (1,1)$, $P_{k}\left[ {%
\begin{array}{c}
\theta \\ 
\phi%
\end{array}%
}\right] (z,\tau )$ is periodic in $z$ with periods $2\pi i\tau $ and $2\pi
i $ with multipliers $\theta $ and $\phi $ respectively. $\ \square $
\end{lemma}

\begin{remark}
\label{Remark_P1_period} Note that the periodicity in $2\pi i$ is determined
by the second argument $\phi $ in contradistinction to the periodicity of
the standard theta series (\ref{Thetaperiod_2pi_i}). Periodicity for $%
(\theta ,\phi )=(1,1)$ is given by (\ref{Pkperiod}).
\end{remark}

We now consider the modular properties. Define the standard left action of
the modular group for $\gamma =\left( 
\begin{array}{cc}
a & b \\ 
c & d%
\end{array}%
\right) \in \Gamma =SL(2,\mathbb{Z})$ on $(z,\tau )\in \mathbb{C}\times 
\mathbb{H}$ with 
\begin{equation}
\gamma .(z,\tau )=(\gamma .z,\gamma .\tau )=(\frac{z}{c\tau +d},\frac{a\tau
+b}{c\tau +d}).  \label{mod_z_tau}
\end{equation}%
We also define a \emph{left} action of $\Gamma $ on $(\theta ,\phi )$ 
\begin{equation}
\gamma .\left[ {%
\begin{array}{c}
\theta \\ 
\phi%
\end{array}%
}\right] =\left[ {%
\begin{array}{c}
\theta ^{a}\phi ^{b} \\ 
\theta ^{c}\phi ^{d}%
\end{array}%
}\right] .  \label{mod_uv}
\end{equation}%
Then we obtain:

\begin{proposition}
\label{Proposition_Pk_Modular} For $(\theta ,\phi )\neq (1,1)$ we have 
\begin{equation}
P_{k}\left( \gamma .\left[ {%
\begin{array}{c}
\theta \\ 
\phi%
\end{array}%
}\right] \right) (\gamma .z,\gamma .\tau )=(c\tau +d)^{k}P_{k}\left[ {%
\begin{array}{c}
\theta \\ 
\phi%
\end{array}%
}\right] (z,\tau ).  \label{Pk_mod}
\end{equation}
\end{proposition}

\noindent \textbf{Proof.} Consider the case $k=1$. It is sufficient to
consider the action of the generators $S,T$ of $\Gamma $ where $S.(z,\tau
)=(-\frac{z}{\tau },-\frac{1}{\tau })$ and $T.(z,\tau )=(z,\tau +1)$. Then
for $\gamma =S$, (\ref{Pk_mod}) is given by (\ref{P1mod_S}) whereas for $%
\gamma =T$, the result follows directly from definition (\ref{Pkdef}). It is
straightforward to check the relations $(ST)^{3}=-S^{2}=1$ (using (\ref%
{P1minusz})) so that the result follows for $k=1$. The general case follows
from (\ref{Pkplus1}). $\square $

\begin{remark}
\label{Remark_Pk_mod}

\begin{description}
\item[(i)] (\ref{Pk_mod}) is equivalent to Theorem 4.2 of \cite{DLM1} for
rational $\lambda$ after noting Remark \ref{Remark_Pk} (i) and (\ref{mod_uv}%
).

\item[(ii)] For $\gamma =-I$ one finds 
\begin{equation}
P_{k}\left[ {%
\begin{array}{c}
\theta \\ 
\phi%
\end{array}%
}\right] (z,\tau )=(-1)^{k}P_{k}\left[ {%
\begin{array}{c}
\theta ^{-1} \\ 
\phi ^{-1}%
\end{array}%
}\right] (-z,\tau ).  \label{Pkodd}
\end{equation}
\end{description}
\end{remark}

\bigskip We next introduce twisted Eisenstein series for $n\geq 1$, defined
by 
\begin{eqnarray}
E_{n}\left[ {%
\begin{array}{c}
\theta \\ 
\phi%
\end{array}%
}\right] (\tau ) &=&-\frac{B_{n}(\lambda )}{n!}+\frac{1}{(n-1)!}%
\sum\limits_{r\geq 0}^{\prime }\frac{(r+\lambda )^{n-1}\theta
^{-1}q^{r+\lambda }}{1-\theta ^{-1}q^{r+\lambda }}  \notag \\
&&+\frac{(-1)^{n}}{(n-1)!}\sum\limits_{r\geq 1}\frac{(r-\lambda
)^{n-1}\theta q^{r-\lambda }}{1-\theta q^{r-\lambda }},  \label{Ekuv}
\end{eqnarray}%
where $\sum\limits^{\prime }$ means we omit $r=0$ if $(\theta ,\phi )=(1,1)$
and where $B_{n}(\lambda )$ is the Bernoulli polynomial defined by%
\begin{equation}
\frac{q_{z}^{\lambda }}{q_{z}-1}=\frac{1}{z}+\sum\limits_{n\geq 1}\frac{%
B_{n}(\lambda )}{n!}z^{n-1}.  \label{Bernoulli}
\end{equation}%
In particular, we note that $B_{1}(\lambda )=\lambda -\frac{1}{2}$.

\begin{remark}
\label{Remark_Ek}\ (i) (\ref{Ekuv}) was introduced in \cite{DLM1} for
rational $\lambda $ where it was denoted by $Q_{n}(\phi ,\theta ^{-1},\tau )$%
.

(ii) $E_{n}\left[ 
\begin{array}{c}
1 \\ 
1%
\end{array}%
\right] (\tau )=E_{n}(\tau )$, the standard Eisenstein series for even $%
n\geq 2$, whereas $E_{n}\left[ 
\begin{array}{c}
1 \\ 
1%
\end{array}%
\right] (\tau )=-B_{1}(0)\delta _{n,1}=\frac{1}{2}\delta _{n,1}$ for $n$ odd.
\end{remark}

We may obtain a Laurant expansion analogous to (\ref{Pkdef}).

\begin{proposition}
\label{Proposition_Pk_Laurant} We have 
\begin{equation}
P_{k}\left[ {%
\begin{array}{c}
\theta \\ 
\phi%
\end{array}%
}\right] (z,\tau )=\frac{1}{z^{k}}+(-1)^{k}\sum\limits_{n\geq k}\binom{n-1}{%
k-1}E_{n}\left[ {%
\begin{array}{c}
\theta \\ 
\phi%
\end{array}%
}\right] (\tau )z^{n-k}.  \label{PkEn}
\end{equation}
\end{proposition}

\noindent \textbf{Proof. } Consider (\ref{Pkuv}) for $k=1$: 
\begin{eqnarray*}
P_{1}\left[ {%
\begin{array}{c}
\theta \\ 
\phi%
\end{array}%
}\right] (z,\tau ) &=&-\sum\limits_{r\geq 0}^{\prime }\frac{q_{z}^{r+\lambda
}}{1-\theta ^{-1}q^{r+\lambda }}-\sum\limits_{r\geq 1}\frac{%
q_{z}^{-r+\lambda }}{1-\theta ^{-1}q^{-r+\lambda }} \\
&=&\frac{q_{z}^{\lambda }}{q_{z}-1}-\sum\limits_{r\geq 0}^{\prime
}q_{z}^{r+\lambda }\frac{\theta ^{-1}q^{r+\lambda }}{1-\theta
^{-1}q^{r+\lambda }} \\
&&+\sum\limits_{r\geq 1}q_{z}^{-r+\lambda }\frac{\theta q^{r-\lambda }}{%
1-\theta q^{r-\lambda }} \\
&=&\frac{1}{z}-\sum\limits_{n\geq 1}E_{n}\left[ {%
\begin{array}{c}
\theta \\ 
\phi%
\end{array}%
}\right] (\tau )z^{n-1},
\end{eqnarray*}%
from (\ref{Ekuv}) and (\ref{Bernoulli}). The general result then follows
from (\ref{Pkplus1}). $\square $

\begin{remark}
\label{Remark_Pk_laurant}\ For $(\theta ,\phi )=(1,1)$ we have $P_{k}\left[ {%
\begin{array}{c}
1 \\ 
1%
\end{array}%
}\right] (z,\tau )=\frac{1}{2}\delta _{k,1}+P_{k}(z,\tau )$ for $k\geq 1$.
\end{remark}

We also find

\begin{proposition}
\label{Proposition_Ek_sum} For $\phi \neq 1$ then%
\begin{equation}
E_{k}\left[ {%
\begin{array}{c}
\theta \\ 
\phi%
\end{array}%
}\right] (\tau )=\frac{1}{(2\pi i)^{k}}\sum_{m\in \mathbb{Z}}\theta ^{m}%
\left[ \sum_{\substack{ n\in \mathbb{Z}  \\ (m,n)\neq (0,0)}}\frac{\phi ^{n}%
}{(m\tau +n)^{k}}\right] ,  \label{Eksum1}
\end{equation}%
whereas for $\theta \neq 1$ 
\begin{equation}
E_{k}\left[ {%
\begin{array}{c}
\theta \\ 
\phi%
\end{array}%
}\right] (\tau )=\frac{1}{(2\pi i)^{k}}\sum_{n\in \mathbb{Z}}\phi ^{n}\left[
\sum_{\substack{ m\in \mathbb{Z}  \\ (m,n)\neq (0,0)}}\frac{\theta ^{m}}{%
(m\tau +n)^{k}}\right] .  \label{Eksum2}
\end{equation}
\end{proposition}

\noindent \textbf{Proof. }Expand the sum of (\ref{P1_double1}) for $\phi
\neq 1$ for $k=1~$to find%
\begin{equation*}
P_{1}\left[ {%
\begin{array}{c}
\theta \\ 
\phi%
\end{array}%
}\right] (z,\tau )=\frac{1}{z}-\frac{1}{2\pi i}\sum_{m\in \mathbb{Z}}\theta
^{m}\left[ \sum_{\substack{ n\in \mathbb{Z}  \\ (m,n)\neq (0,0)}}%
\sum\limits_{r\geq 1}(\frac{z}{2\pi i})^{r-1}\frac{\phi ^{n}}{(m\tau +n)^{r}}%
\right] .
\end{equation*}%
Comparing with (\ref{PkEn}) then (\ref{Eksum1}) follows. (\ref{Eksum2})
similarly holds. $\square $

\begin{remark}
\label{Remark_Eksum}\ When both $\theta \neq 1$ and $\phi \neq 1$ then (\ref%
{Eksum1}) and (\ref{Eksum2}) are equal. For $k\geq 3$, they are absolutely
convergent and equal for all $(\theta ,\phi )$. For $k\geq 3$, and $(\theta
,\phi )=(1,1)$ we obtain the standard Eisenstein series (\ref{Eisen}).
\end{remark}

From Proposition \ref{Proposition_Pk_Modular} it immediately follows that

\begin{proposition}
\label{Proposition_Ek_modular} For $(\theta ,\phi )\neq (1,1)$, $E_{k}\left[ 
{%
\begin{array}{c}
\theta \\ 
\phi%
\end{array}%
}\right] $ is a modular form of weight $k$ where 
\begin{equation}
E_{k}\left( \gamma .\left[ {%
\begin{array}{c}
\theta \\ 
\phi%
\end{array}%
}\right] \right) (\gamma .\tau )=(c\tau +d)^{k}E_{k}\left[ {%
\begin{array}{c}
\theta \\ 
\phi%
\end{array}%
}\right] (\tau ).\quad \square  \label{Ekmodular}
\end{equation}
\end{proposition}

\begin{remark}
\label{Remark_Ek_usual}\ This is equivalent to Theorem 4.6 of \cite{DLM1}
for rational $\lambda $. (\ref{Ekmodular}) also holds for $(\theta ,\phi
)=(1,1)$ for $k\geq 3$, whereas $E_{2}$ is quasi-modular.
\end{remark}

\medskip It is useful to note the analytic expansions: 
\begin{eqnarray}
P_{1}\left[ {%
\begin{array}{c}
\theta \\ 
\phi%
\end{array}%
}\right] (z_{1}-z_{2},\tau ) &=&\frac{1}{z_{1}-z_{2}}+\sum_{k,l\geq 1}C\left[
{%
\begin{array}{c}
\theta \\ 
\phi%
\end{array}%
}\right] (k,l)z_{1}^{k-1}z_{2}^{l-1},  \label{P1_C_expansion} \\
P_{1}\left[ {%
\begin{array}{c}
\theta \\ 
\phi%
\end{array}%
}\right] (z+z_{1}-z_{2},\tau ) &=&\sum_{k,l\geq 1}D\left[ {%
\begin{array}{c}
\theta \\ 
\phi%
\end{array}%
}\right] (k,l,z)z_{1}^{k-1}z_{2}^{l-1},  \label{P1_D_expansion}
\end{eqnarray}%
where for $k,l\geq 1$ we define 
\begin{eqnarray}
C\left[ {%
\begin{array}{c}
\theta \\ 
\phi%
\end{array}%
}\right] (k,l,\tau ) &=&(-1)^{l}\binom{k+l-2}{k-1}E_{k+l-1}\left[ {%
\begin{array}{c}
\theta \\ 
\phi%
\end{array}%
}\right] (\tau ),  \label{Ckldef} \\
D\left[ {%
\begin{array}{c}
\theta \\ 
\phi%
\end{array}%
}\right] (k,l,\tau ,z) &=&(-1)^{k+1}\binom{k+l-2}{k-1}P_{k+l-1}\left[ {%
\begin{array}{c}
\theta \\ 
\phi%
\end{array}%
}\right] (\tau ,z).  \label{Dkldef}
\end{eqnarray}%
We also note that (\ref{Pkodd}) implies 
\begin{eqnarray}
C\left[ {%
\begin{array}{c}
\theta \\ 
\phi%
\end{array}%
}\right] (k,l,\tau ) &=&-C\left[ {%
\begin{array}{c}
\theta ^{-1} \\ 
\phi ^{-1}%
\end{array}%
}\right] (l,k,\tau ),  \label{Cklodd} \\
D\left[ {%
\begin{array}{c}
\theta \\ 
\phi%
\end{array}%
}\right] (k,l,\tau ,z) &=&-D\left[ {%
\begin{array}{c}
\theta ^{-1} \\ 
\phi ^{-1}%
\end{array}%
}\right] (l,k,\tau ,-z).  \label{Dklodd}
\end{eqnarray}

Finally, we may also express the twisted Weierstrass functions in terms of
theta series and the prime form as follows:

\begin{proposition}
\label{Proposition_Weier_theta}For $(\theta ,\phi )\neq (1,1)$ with $\theta
=\exp (-2\pi i\mu )$ and $\phi =\exp (2\pi i\lambda )$ then 
\begin{equation}
P_{1}\left[ {%
\begin{array}{c}
\theta \\ 
\phi%
\end{array}%
}\right] (z,\tau )=\frac{\vartheta \left[ {%
\begin{array}{c}
\lambda +\frac{1}{2} \\ 
\mu +\frac{1}{2}%
\end{array}%
}\right] (z,\tau )}{\vartheta \left[ {%
\begin{array}{c}
\lambda +\frac{1}{2} \\ 
\mu +\frac{1}{2}%
\end{array}%
}\right] (0,\tau )}\frac{1}{K(z,\tau )},  \label{P1uv_theta}
\end{equation}%
whereas 
\begin{equation}
P_{1}\left[ {%
\begin{array}{c}
1 \\ 
1%
\end{array}%
}\right] (z,\tau )=\frac{\frac{d}{dz}\vartheta \left[ {%
\begin{array}{c}
\frac{1}{2} \\ 
\frac{1}{2}%
\end{array}%
}\right] (z,\tau )}{\frac{d}{dz}\vartheta \left[ {%
\begin{array}{c}
\frac{1}{2} \\ 
\frac{1}{2}%
\end{array}%
}\right] (0,\tau )}\frac{1}{K(z,\tau )}.  \label{P1theta}
\end{equation}
\end{proposition}

\noindent \textbf{Proof.} For $(\theta ,\phi )\neq (1,1)$ the result follows
by comparing the periodicity and pole structure of each expression using (%
\ref{Thetaperiod_2pi_i}) and (\ref{Thetaperiod_2Pi_i_tau}). For $(\theta
,\phi )=(1,1)$ the result follows from (\ref{P1_dlogKdz}) and (\ref{Ktheta}%
). $\square $


\section{\label{Section_vosas}$n$-Point Functions for $\mathbb{R}$-Graded
Vertex Operator Superalgebras}

\subsection{Introduction to Vertex Operator Superalgebras}

We discuss some aspects of Vertex Operator Superalgebra (VOSA) theory to
establish context and notation. For more details see \cite{B}, \cite{FHL}, 
\cite{FLM}, \cite{Ka}, \cite{MN1}. Let $V$ be a superspace i.e. a complex
vector space $V=V_{\bar{0}}\oplus V_{\bar{1}}=\oplus _{\alpha }V_{\alpha }$
with index label $\alpha $ in $\mathbb{Z}/2\mathbb{Z}$ so that each $a\in V$
has a parity (fermion number) $p(a)\in \mathbb{Z}/2\mathbb{Z}$.

An $\mathbb{R}$-graded Vertex Operator Superalgebra (VOSA) is a quadruple $%
(V,Y,\mathbf{1},\omega )$ as follows: $V$ is a superspace with a (countable) 
$\mathbb{R}$-grading where 
\begin{equation*}
V=\oplus _{r\geq r_{0}}V_{r}
\end{equation*}%
for some $r_{0}$ and with parity decomposition $V_{r}=V_{\bar{0},r}\oplus V_{%
\bar{1},r}$. $\mathbf{1}\in V_{\bar{0},0}$ is the vacuum vector and $\omega
\in V_{\bar{0},2}$ the conformal vector with properties described below. $Y$
is a linear map $Y:V\rightarrow (\mathrm{End}V)[[z,z^{-1}]]$, for formal
variable $z$, so that for any vector (state) $a\in V$ 
\begin{equation}
Y(a,z)=\sum_{n\in \mathbb{Z}}a(n)z^{-n-1}.  \label{Ydefn}
\end{equation}%
The component operators (modes) $a(n)\in \mathrm{End}V$ are such that $a(n)%
\mathbf{1}=\delta _{n,-1}a$ for $n\geq -1$ and 
\begin{equation}
a(n)V_{\alpha }\subset V_{\alpha +p(a)},  \label{a(n)parity}
\end{equation}%
for $a\ $of parity $p(a)$.

The vertex operators satisfy the locality property for all $a,b\in V$%
\begin{equation}
(x-y)^{N}[Y(a,x),Y(b,y)]=0,  \label{locality}
\end{equation}%
for $N\gg 0$, where the commutator is defined in the graded sense, i.e. 
\begin{equation*}
\lbrack Y(a,x),Y(b,y)]=Y(a,x)Y(b,y)-(1)^{p(a)p(b)}Y(b,y)Y(a,x).
\end{equation*}

The vertex operator for the vacuum is $Y(\mathbf{1},z)=Id_{V}$, whereas that
for $\omega $ is 
\begin{equation}
Y(\omega ,z)=\sum_{n\in \mathbb{Z}}L(n)z^{-n-2},  \label{Yomega}
\end{equation}%
where $L(n)$ forms a Virasoro algebra for central charge $c$%
\begin{equation}
\lbrack L(m),L(n)]=(m-n)L(m+n)+\frac{c}{12}(m^{3}-m)\delta _{m,-n}.
\label{Virasoro}
\end{equation}%
$L(-1)$ satisfies the translation property 
\begin{equation}
Y(L(-1)a,z)=\frac{d}{dz}Y(a,z).  \label{YL(-1)}
\end{equation}%
$L(0)$ describes the $\mathbb{R}$-grading with $L(0)a=wt(a)a$ for weight $%
wt(a)\in \mathbb{R}$ and 
\begin{equation}
V_{r}=\{a\in V|wt(a)=r\}.  \label{Vdecomp}
\end{equation}

We quote the standard commutator property of VOSAs e.g. \cite{Ka}, \cite{FHL}%
, \cite{MN1} 
\begin{equation}
\lbrack a(m),Y(b,z)]=\sum\nolimits_{j\geq 0}\binom{m}{j}Y(a(j)b,z)z^{m-j}.
\label{aYcomm}
\end{equation}%
Taking $a=\omega $ this implies for $b$ of weight $wt(b)$ that 
\begin{equation}
\lbrack L(0),b(n)]=(wt(b)-n-1)b(n),  \label{L0b}
\end{equation}%
so that 
\begin{equation}
b(n)V_{r}\subset V_{r+wt(b)-n-1}.  \label{bnVr}
\end{equation}%
In particular, we define for $a$ of weight $wt(a)$ the zero mode 
\begin{equation}
o(a)=\left\{ 
\begin{array}{c}
a(wt(a)-1),\text{ \ for }wt(a)\in \mathbb{Z} \\ 
0\text{, \ otherwise,}%
\end{array}%
\right.  \label{o(v)}
\end{equation}%
which is then extended by linearity to all $a\in V$.

\subsection{Torus $n$-point Functions}

In this section we will develop explicit formulas for the $n$-point
functions for $\mathbb{R}$-graded VOSA modules at genus one \cite{Z, DLM1,
MT1, DZ1}. Let $(V,Y,\mathbf{1},\omega )$ be an $\mathbb{R}$-graded VOSA. In
order to consider modular-invariance of $n$-point functions at genus 1, Zhu
introduced in ref. \cite{Z} a second "square-bracket" VOA $(V,Y[,],\mathbf{1}%
,\tilde{\omega})$ associated to a given VOA $(V,Y(,),\mathbf{1},\omega )$.
We review some aspects of that construction here. The new square bracket
vertex operators are defined by a change of co-ordinates, namely 
\begin{equation}
Y[v,z]=\sum_{n\in \mathbb{Z}}v[n]z^{-n-1}=Y(q_{z}^{L(0)}v,q_{z}-1),
\label{Ysquare}
\end{equation}%
with $q_{z}=\exp (z)$, while the new conformal vector is $\tilde{\omega}%
=\omega -\frac{c}{24}\mathbf{1}$. 
For $v$ of $L(0)$ weight $wt(v)\in \mathbb{R}$ and $m\geq 0$, 
\begin{eqnarray}
v[m] &=&m!\sum\limits_{i\geq m}c(wt(v),i,m)v(i),  \label{square1} \\
\sum\limits_{m=0}^{i}c(wt(v),i,m)x^{m} &=&\binom{wt(v)-1+x}{i}.
\label{square2}
\end{eqnarray}%
In particular we note that $v[0]=\sum\limits_{i\geq 0}\binom{wt(v)-1}{i}v(i)$%
. 
\qquad

\medskip We now define the torus $n$-point functions. Following (\ref%
{a(n)parity}) we let $\sigma \in \mathrm{Aut}(V)$ denote the parity (fermion
number) automorphism 
\begin{equation}
\sigma a=(-1)^{p(a)}a.  \label{sigma}
\end{equation}%
Let $g\in \mathrm{Aut}(V)$ denote any other automorphism which commutes with 
$\sigma $. Let $M$ be a $V$-module with vertex operators $Y_{M}$. We assume
that $M$ is stable under both $\sigma $ and $g$ i.e. $\sigma $ and $g$ act
on $M$ \cite{DZ1}. The $n$\emph{-point function} on $M$ for states $%
v_{1},\ldots ,v_{n}\in V$ and $g\in \mathrm{Aut}(V)$ is defined by\footnote{%
This $n$-point function would be denoted by $T((v_{1},q_{1}),\ldots
,(v_{n},q_{n}),(1,g),q)$ in the notation of \cite{DLM1} and \cite{DZ1}.} 
\begin{gather}
F_{M}(g;v_{1},\ldots v_{n};\tau )=F_{M}(g;(v_{1},z_{1}),\ldots
,(v_{n},z_{n});\tau )  \notag \\
=\mathrm{STr}_{M}\left( g\;Y_{M}(q_{1}^{L(0)}v_{1},q_{1})\ldots
Y_{M}(q_{n}^{L(0)}v_{n},q_{n})q^{L(0)-c/24}\right) ,  \label{npointfunction}
\end{gather}%
$q=\exp (2\pi i\tau )$, $q_{i}=\exp (z_{i})$, $1\leq i\leq n$, for auxiliary
variables $z_{1},...,z_{n}$ and where $\mathrm{STr}_{M}$ denotes the \emph{%
supertrace} defined by%
\begin{equation}
\mathrm{STr}_{M}(X)=Tr_{M}(\sigma X)=Tr_{M_{\bar{0}}}(X)-Tr_{M_{\bar{1}}}(X).
\label{Supertrace}
\end{equation}%
In Appendix A we describe some basic properties of the supertrace. Taking $%
g=1$ and all $v_{i}=\mathbf{1}$ in (\ref{npointfunction}) yields the \emph{%
partition function} which we denote by 
\begin{equation}
Z_{M}(\tau )=F_{M}(1;\tau )=\mathrm{STr}_{M}\left( q^{L(0)-c/24}\right) .
\label{ZM}
\end{equation}%
We also denote the orbifold partition function for general $g$ by 
\begin{equation}
Z_{M}(g,\tau )=F_{M}(g;\tau )=\mathrm{STr}_{M}\left( gq^{L(0)-c/24}\right) .
\label{ZMg}
\end{equation}

For $g=1$ (\ref{npointfunction}) is defined by Zhu for a $\mathbb{Z}$-graded
VOA \cite{Z}. For $g$ of finite order, it is considered for $\mathbb{Z}$%
-graded VOAs in ref. \cite{DLM1}, $\frac{1}{2}\mathbb{Z}$-graded VOSAs in
ref. \cite{DZ1} and $\mathbb{Z}$-graded VOSAs in ref. \cite{DZ2}. Here we
generalize these results to an $\mathbb{R}$-graded VOSA for arbitrary $g$
commuting with $\sigma $.

For $n=1$ in (\ref{npointfunction}) we obtain the \emph{1-point function}
denoted by

\begin{equation}
Z_{M}(g,v_{1},\tau )=F_{M}(g;(v_{1},z_{1});\tau )=\mathrm{STr}%
_{M}(go(v_{1})q^{L(0)-c/24})  \label{ZM1pt}
\end{equation}%
where $o(v_{1})$ is the zero mode (\ref{o(v)}) and is independent of $z_{1}$%
. We note the following useful result relating any $n$-point function to a
1-point function:

\begin{lemma}
\label{lemma_npt_1pt} For states $v_{1},v_{2},\ldots ,v_{n}$ as above we have%
\begin{eqnarray}
&&F_{M}(g;(v_{1},z_{1}),\ldots ,(v_{n},z_{n});\tau )  \notag \\
&=&Z_{M}(g,Y[v_{1},z_{1n}].Y[v_{2},z_{2n}]\ldots
Y[v_{n-1},z_{n-1n}].v_{n},\tau )  \label{Fnziminuszn} \\
&=&Z_{M}(g,Y[v_{1},z_{1}].Y[v_{2},z_{2}]\ldots Y[v_{n},z_{n}].\mathbf{1}%
,\tau ),  \label{Fnz1zn}
\end{eqnarray}%
where $z_{ij}=z_{i}-z_{j}$.
\end{lemma}

\noindent \textbf{Proof. }The proof follows Lemma 1 of ref. \cite{MT1}. $%
\square $

\medskip Every $n$-point function enjoys the following permutation and
periodicity properties \cite{Z}, \cite{MT1}:

\begin{lemma}
\label{Lemma_Perm_Period} Consider the $n$-point function $F_{M}$ for states 
$v_{1},v_{2},\ldots ,v_{n}$, as above, where each $v_{i}$ is of weight $%
wt(v_{i})$, parity $p(v_{i})$ and is a $g$-eigenvector for eigenvalue $%
\theta _{i}^{-1}$.

\begin{description}
\item[(i)] If $p(v_{1})+\ldots +p(v_{n})$ is odd then $F_{M}=0$.

\item[(ii)] Permuting adjacent vectors, 
\begin{eqnarray*}
&&F_{M}(g;(v_{1},z_{1}),\ldots ,(v_{k},z_{k}),(v_{k+1},z_{k+1}),\ldots
,(v_{n},z_{n});\tau ) \\
&=&(-1)^{p(v_{k})p(v_{k+1})}F_{M}(g;(v_{1},z_{1}),\ldots
,(v_{k+1},z_{k+1}),(v_{k},z_{k}),\ldots ,(v_{n},z_{n});\tau ).
\end{eqnarray*}

\item[(iii)] $F_{M}$ is a function of $z_{ij}=$ $z_{i}-z_{j}$ and is
non-singular at $z_{ij}$ $\neq 0$ for all $i\neq j$.

\item[(iv)] $F_{M}$ is periodic in $z_{i}$ with period $2\pi i$ and
multiplier $\phi _{i}=\exp (2\pi iwt(v_{i}))$.

\item[(v)] $F_{M}$ is periodic in $z_{i}$ with period $2\pi i\tau $ and
multiplier $\theta _{i}$.
\end{description}
\end{lemma}

\noindent \textbf{Proof}. (i) This follows from definition (\ref{Supertrace}%
).

(ii) Apply locality (\ref{locality}).

(iii) $F_{M}$ is a function of $z_{ij}$ from (\ref{Fnziminuszn}). Suppose $%
F_{M}$ is singular at $z_{n}=y$ for some $y\neq $ $z_{j}$ for all $%
j=1,\ldots ,n-1$. We may assume that $z_{0}=0$ by redefining $z_{i}$ to be $%
z_{i}-z_{0}$ for all $i$. But from (\ref{Fnz1zn}), $Y[v_{n},z_{n}].\mathbf{1}%
|_{z_{n}=0}=v_{n}$ is non-singular at $z_{nj}$ $\neq $ $0$. Applying (ii)
the result follows for all $z_{ij}$.

(iv) This follows directly from the definition (\ref{npointfunction}).

(v) Using (iii) we consider periodicity of $z_{n}$ wlog. Under $%
z_{n}\rightarrow z_{n}+2\pi i\tau \,$ $\,$we have $F_{M}\rightarrow \hat{F}%
_{M}$ where 
\begin{eqnarray*}
\hat{F}_{M} &=&q^{-c/24}\mathrm{STr}_{M}(gY(q_{1}^{L(0)}v_{1},q_{1})\ldots
Y(q^{L(0)}q_{n}^{L(0)}v_{n},qq_{n})q^{L(0)}) \\
&=&q^{-c/24}\mathrm{STr}_{M}(gY(q_{1}^{L(0)}v_{1},q_{1})\ldots
q^{L(0)}Y(q_{n}^{L(0)}v_{n},q_{n})),
\end{eqnarray*}%
using $q^{L(0)}Y(b,z)q^{-L(0)}=Y(q^{L(0)}b,qz)$ (which follows from (\ref%
{L0b})). But 
\begin{eqnarray*}
&&\mathrm{STr}_{M}(gY(q_{1}^{L(0)}v_{1},q_{1})\ldots
q^{L(0)}Y(q_{n}^{L(0)}v_{n},q_{n})) \\
&=&(-1)^{p(v_{n})}\mathrm{STr}%
_{M}(Y(q_{n}^{L(0)}v_{n},q_{n})gY(q_{1}^{L(0)}v_{1},q_{1})\ldots
Y(q_{n-1}^{L(0)}v_{n-1},q_{n-1})q^{L(0)}) \\
&=&\theta _{n}(-1)^{p(v_{n})}\mathrm{STr}%
_{M}(gY(q_{n}^{L(0)}v_{n},q_{n})Y(q_{1}^{L(0)}v_{1},q_{1})\ldots
Y(q_{n-1}^{L(0)}v_{n-1},q_{n-1})q^{L(0)}) \\
&=&\theta _{n}\mathrm{STr}_{M}(gY(q_{1}^{L(0)}v_{1},q_{1})\ldots
Y(q_{n-1}^{L(0)}v_{n-1},q_{n-1})Y(q_{n}^{L(0)}v_{n},q_{n})q^{L(0)}),
\end{eqnarray*}%
using $g^{-1}Y(v_{n},q_{n})g=Y(g^{-1}v_{n},q_{n})=\theta _{n}Y(v_{n},q_{n})$
and applying (ii) repeatedly. Thus $\hat{F}_{M}=\theta _{n}F_{M}$. \ $%
\square $

\subsection{Zhu Recursion Formulas for $n$-Point Functions}

\label{apb}

We now prove a generalization of Zhu's $n$-point function recursion formula 
\cite{Z} for the $n$-point function (\ref{npointfunction}) for an $\mathbb{R}
$-graded VOSA. We begin with the following Lemma which follows directly from
(\ref{aYcomm}):

\begin{lemma}
\label{lemmaacomm} Suppose that $u\in V$ is homogeneous of weight $wt(u)\in 
\mathbb{R}$. Then for $k\in \mathbb{Z}$ and $v\in V$ we have 
\begin{equation}
\left[ u(k),Y(q_{z}^{L(0)}v,q_{z})\right] =q_{z}^{k-wt(u)+1}\sum\limits_{i%
\geq 0}\binom{k}{i}Y(q_{z}^{L(0)}u(i)v,q_{z}).\square  \label{akcomm}
\end{equation}
\end{lemma}

\begin{corollary}
\label{corro(a)comm}Suppose that $u\in V$ is homogeneous of integer weight $%
wt(u)\in \mathbb{Z}$. Then 
\begin{equation}
\left[ o(u),Y(q_{z}^{L(0)}v,q_{z})\right] =Y(q_{z}^{L(0)}u[0]v,q_{z}).
\label{o(u)comm}
\end{equation}
\end{corollary}

Similarly to Zhu's Proposition 4.3.1 (op.cit.) we find

\begin{proposition}
\label{Propa[0]comm} Suppose that $v\in V$ is homogeneous of integer weight $%
wt(v)\in \mathbb{Z}$. Then for $v_{1},\ldots ,v_{n}\in V$, we have 
\begin{equation}
\sum\limits_{r=1}^{n}p(v,v_{1}v_{2}\ldots v_{r-1})F_{M}(g;v_{1};\ldots
;v[0]v_{r};\ldots v_{n};\tau )=0,  \label{v[0]comm}
\end{equation}%
with $p(v,v_{1}v_{2}\ldots v_{r-1})$ of (\ref{parityAB}) in Appendix A. $%
\square $
\end{proposition}

Let $v$ be homogeneous of weight $wt(v)\in \mathbb{R}$ and define $\phi \in
U(1)$ by 
\begin{equation}
\phi =\exp (2\pi iwt(v)).  \label{phi}
\end{equation}%
We also take $v$ to be an eigenfunction under $g$ with 
\begin{equation}
gv=\theta ^{-1}v  \label{theta}
\end{equation}%
for some $\theta \in U(1)$ so that 
\begin{equation}
g^{-1}v(k)g=\theta v(k).  \label{gv(k)}
\end{equation}%
Then we obtain the following generalization of Zhu's Proposition 4.3.2 \cite%
{Z} for the $n$-point function:\ 

\begin{theorem}
\label{Theorem_npt_rec}Let $v,\theta $ and $\phi $\ be as as above. Then for
any $v_{1},\ldots v_{n}\in V$ we have 
\begin{eqnarray}
&&F_{M}(g;v,v_{1},\ldots v_{n};\tau )  \notag \\
&=&\delta _{\theta ,1}\delta _{\phi ,1}\mathrm{STr}_{M}\left(
go(v)Y_{M}(q_{1}^{L(0)}v_{1},q_{1})\ldots
Y_{M}(q_{n}^{L(0)}v_{n},q_{n})q^{L(0)-c/24}\right)  \notag \\
&&+\sum\limits_{r=1}^{n}\sum\limits_{m\geq 0}p(v,v_{1}v_{2}\ldots
v_{r-1})P_{m+1}\left[ 
\begin{array}{c}
\theta \\ 
\phi%
\end{array}%
\right] (z-z_{r},\tau )\cdot  \notag \\
&&F_{M}(g;v_{1},\ldots ,v[m]v_{r},\ldots ,v_{n};\tau ).  \label{nptrec}
\end{eqnarray}%
(The twisted Weierstrass function is defined in (\ref{Pkuv})).
\end{theorem}

\noindent \textbf{Proof.} We have%
\begin{eqnarray*}
&&q^{c/24}F_{M}(g;v,v_{1},\ldots v_{n};\tau ) \\
&=&\sum_{k\in \mathbb{Z}}q_{z}^{-k-1+wt(v)}\mathrm{STr}_{M}\left(
g\;v(k)Y_{M}(q_{1}^{L(0)}v_{1},q_{1})\ldots
Y_{M}(q_{n}^{L(0)}v_{n},q_{n})q^{L(0)}\right) .
\end{eqnarray*}

Thus we consider 
\begin{eqnarray*}
&&\mathrm{STr}_{M}\left( g\;v(k)Y_{M}(q_{1}^{L(0)}v_{1},q_{1})\ldots
Y_{M}(q_{n}^{L(0)}v_{n},q_{n})q^{L(0)}\right) \\
&=&\mathrm{STr}_{M}\left( g\;[v(k),Y_{M}(q_{1}^{L(0)}v_{1},q_{1})\ldots
Y_{M}(q_{n}^{L(0)}v_{n},q_{n})]q^{L(0)}\right) \\
&&+p(v,v_{1}\cdots v_{n})\mathrm{STr}_{M}\left(
g\;Y_{M}(q_{1}^{L(0)}v_{1},q_{1})\ldots
Y_{M}(q_{n}^{L(0)}v_{n},q_{n})v(k)q^{L(0)}\right) \\
&=&\sum_{r=1}^{n}\sum_{i\geq 0}p(v,v_{1}\ldots v_{r-1})\binom{k}{i}%
q_{r}^{k+1-wt(v)}. \\
&&\mathrm{STr}_{M}\left( g\;Y_{M}(q_{1}^{L(0)}v_{1},q_{1})\ldots
Y_{M}(q_{r}^{L(0)}v(i)v_{r},q_{r})\ldots
Y_{M}(q_{n}^{L(0)}v_{n},q_{n})q^{L(0)}\right) \\
&&+\theta q^{k+1-wt(v)}\mathrm{STr}_{M}\left(
g\;v(k)Y_{M}(q_{1}^{L(0)}v_{1},q_{1})\ldots
Y_{M}(q_{n}^{L(0)}v_{n},q_{n})q^{L(0)}\right) ,
\end{eqnarray*}%
applying (\ref{L0b}), (\ref{akcomm}), (\ref{gv(k)}), (\ref{AB1toBncom}) and
Lemma \ref{lemmacommk} of Appendix A. Thus 
\begin{eqnarray*}
&&\mathrm{STr}_{M}\left( g\;v(k)Y_{M}(q_{1}^{L(0)}v_{1},q_{1})\ldots
Y_{M}(q_{n}^{L(0)}v_{n},q_{n})q^{L(0)}\right) \\
&=&\frac{1}{1-\theta q^{k+1-wt(v)}}\sum_{r=1}^{n}\sum_{i\geq
0}p(v,v_{1}\ldots v_{r-1})\binom{k}{i}q_{r}^{k-wt(v)+1}. \\
&&\mathrm{STr}_{M}\left( gY(q_{1}^{L(0)}v_{1},q_{1})\ldots
Y(q_{r}^{L(0)}v(i)v_{r},q_{r})\ldots
Y(q_{n}^{L(0)}v_{n},q_{n})q^{L(0)}\right) ,
\end{eqnarray*}%
provided $(\theta ,\phi ,k)\neq (1,1,-1+wt(v))$. This implies $%
F_{M}(g;v,v_{1},\ldots v_{n})$ is given by 
\begin{eqnarray*}
&&\delta _{\theta ,1}\delta _{\phi ,1}\mathrm{STr}_{M}\left(
go(v)Y_{M}(q_{1}^{L(0)}v_{1},q_{1})\ldots
Y_{M}(q_{n}^{L(0)}v_{n},q_{n})q^{L(0)-c/24}\right) \\
&&+\sum_{r=1}^{n}p(v,v_{1}\ldots v_{r-1})\sum_{k\in \mathbb{Z}}^{\prime }%
\frac{(\frac{q_{r}}{q_{z}})^{k+1-wt(v)}}{1-\theta q^{k+1-wt(v)}}%
.F_{M}(g;v_{1},\ldots \sum_{i\geq 0}\binom{k}{i}v(i)v_{r},\ldots ,v_{n}),
\end{eqnarray*}%
where the prime denotes the omission of $k=-1+wt(v)$ if $(\theta ,\phi
)=(1,1)$ and recalling (\ref{o(v)}). Now from (\ref{square1}) and (\ref%
{square2}) we find 
\begin{equation*}
\sum_{i\geq 0}\binom{k}{i}v(i)=\sum\limits_{m\geq 0}\frac{(k+1-wt(v))^{m}}{m!%
}v[m].
\end{equation*}%
The sum over $k$ can then be computed in terms of a twisted Weierstrass
function (\ref{Pkuv}) for $\lambda =wt(v)(\mathrm{mod}\mathbb{Z})$ as
follows: 
\begin{eqnarray*}
&&\frac{1}{m!}\sum_{k\in \mathbb{Z}}^{\prime }\frac{(k+1-wt(v))^{m}(\frac{%
q_{r}}{q_{z}})^{k+1-wt(v)}}{1-\theta q^{k+1-wt(v)}} \\
&=&(-1)^{m+1}P_{m+1}\left[ 
\begin{array}{c}
\theta ^{-1} \\ 
\phi ^{-1}%
\end{array}%
\right] (z_{r}-z,\tau )-\frac{1}{2}\delta _{\theta ,1}\delta _{\phi
,1}\delta _{m,0} \\
&=&P_{m+1}\left[ 
\begin{array}{c}
\theta \\ 
\phi%
\end{array}%
\right] (z-z_{r},\tau )-\frac{1}{2}\delta _{\theta ,1}\delta _{\phi
,1}\delta _{m,0},
\end{eqnarray*}%
using (\ref{Pkodd}). Thus we find $F_{M}(g;v,v_{1},\ldots v_{n},\tau )$ is
given by 
\begin{eqnarray*}
&&\delta _{\theta ,1}\delta _{\phi ,1}\mathrm{STr}_{M}\left(
go(v)Y_{M}(q_{1}^{L(0)}v_{1},q_{1})\ldots
Y_{M}(q_{n}^{L(0)}v_{n},q_{n})q^{L(0)-c/24}\right) \\
&&+\sum\limits_{r=1}^{n}\sum\limits_{m\geq 0}p(v,v_{1}v_{2}\ldots
v_{r-1})P_{m+1}\left[ 
\begin{array}{c}
\theta \\ 
\phi%
\end{array}%
\right] (z-z_{r},\tau )F_{M}(g;v_{1},\ldots ,v[m]v_{r},\ldots ,v_{n};\tau )
\\
&&-\frac{1}{2}\delta _{\theta ,1}\delta _{\phi
,1}\sum_{r=1}^{n}p(v,v_{1}\ldots v_{r-1})F_{M}(g;v_{1},\ldots
,v[0]v_{r},\ldots ,v_{n};\tau ).
\end{eqnarray*}%
Finally, it follows from (\ref{v[0]comm}) that the last sum is zero and
hence (\ref{nptrec}) obtains. $\square $

\begin{remark}
\label{Remark_Cgrading} (i) Note that it is necessary for the $V$-grading to
be real in order in order for $P_{k}\left[ 
\begin{array}{c}
\theta \\ 
\phi%
\end{array}%
\right] $ to converge. Thus, VOSAs with $\mathbb{C}$-grading such as those
discussed in \cite{DM} have divergent torus $n$-point functions.

(ii) From Lemma \ref{Lemma_periodicity} it follows that $F_{M}$ is periodic
in $z$ with periods $2\pi i\tau $ and $2\pi i$ with multipliers $\theta $
and $\phi $ respectively in agreement with Lemma \ref{Lemma_Perm_Period}.
\end{remark}

Other standard recursion formulas can be similarly generalized. Thus

\begin{proposition}
\label{propapnpt} With notation as above, for any states $v_{1},\ldots
v_{n}\in V$, and for $p\geq 1$ we have: 
\begin{eqnarray}
&&F_{M}(g;v[-p].v_{1},\ldots v_{n};\tau )  \notag \\
&=&\delta _{\theta ,1}\delta _{\phi ,1}\delta _{p,1}\;\mathrm{STr}%
_{M}(go(v)Y(q_{1}^{L(0)}v_{1},q_{1})\ldots
Y(q_{n}^{L(0)}v_{n},q_{n})q^{L(0)-c/24})  \notag \\
&&+\sum\limits_{m\geq 0}(-1)^{m+1}\binom{m+p-1}{m}E_{m+p}\left[ 
\begin{array}{c}
\theta \\ 
\phi%
\end{array}%
\right] (\tau )F_{M}(g;v[m]v_{1},\ldots v_{n};\tau )  \notag \\
&&+\sum_{r=2}^{n}\sum\limits_{m\geq 0}p(v,v_{1}v_{2}\ldots v_{r-1})(-1)^{p+1}%
\binom{m+p-1}{m}P_{m+p}\left[ 
\begin{array}{c}
\theta \\ 
\phi%
\end{array}%
\right] (z_{1r},\tau ).  \notag \\
&&F_{M}(g;v_{1},\ldots v[m]v_{r},\ldots v_{n};\tau ).  \label{npt-recursionp}
\end{eqnarray}
\end{proposition}

\noindent \textbf{Proof.} \ Using (\ref{Fnz1zn}) of Lemma \ref{lemma_npt_1pt}
and associativity of VOSAs (e.g. \cite{FHL}) we have: 
\begin{eqnarray}
&&F_{M}(g;(Y[v,z]v_{1},z_{1}),\ldots (v_{n},z_{n});\tau )  \notag \\
&=&Z_{M}(g,Y[Y[v,z]v_{1},z_{1}]Y[v_{2},z_{2}]\ldots Y[v_{n},z_{n}]\mathbf{1}%
,\tau )  \notag \\
&=&Z_{M}(g,Y[v,z+z_{1}]Y[v_{1},z_{1}]Y[v_{2},z_{2}]\ldots Y[v_{n},z_{n}]%
\mathbf{1},\tau )  \notag \\
&=&F_{M}(g;(v,z+z_{1}),(v_{1},z_{1}),\ldots (v_{n},z_{n});\tau ).
\label{FY[v,z]}
\end{eqnarray}%
Expanding the LHS of (\ref{FY[v,z]}) in $z$ we find that the coefficient of $%
z^{p-1}$ is $F_{M}(v[-p].v_{1},z_{1};\ldots v_{n},z_{n};g;\tau )$. We can
compare this to the expansion of the RHS in $z$ from (\ref{nptrec}) of
Theorem \ref{Theorem_npt_rec}. From (\ref{PkEn}) we find that for $p\geq 1$,
the coefficient of $z^{p-1}$ in $P_{m+1}\left[ 
\begin{array}{c}
\theta \\ 
\phi%
\end{array}%
\right] (z,\tau )$ is $(-1)^{m+1}\binom{m+p-1}{m}E_{m+p}\left[ 
\begin{array}{c}
\theta \\ 
\phi%
\end{array}%
\right] (\tau )$. Furthermore for $r\neq 1$ the coefficient of $z^{p-1}$ in $%
P_{m+1}\left[ 
\begin{array}{c}
\theta \\ 
\phi%
\end{array}%
\right] (z+z_{1r},\tau )$ is given by $(-1)^{p+1}\binom{m+p-1}{m}P_{m+p}%
\left[ 
\begin{array}{c}
\theta \\ 
\phi%
\end{array}%
\right] (z_{1r},\tau )$. Lastly, for $p=1$ the first term of (\ref{nptrec})
also contributes. Thus the stated result follows. $\square $

\section{\label{Section_Heisenberg}Shifted VOSAs and Heisenberg Twisted
Modules}

In this section we discuss the $n$-point functions for an orbifold $g$%
-twisted module for a VOSA where $g$ is a continuous symmetry generated by a
Heisenberg vector. (For definitions and properties of twisted modules we
refer the reader to refs. \cite{Li}, \cite{DLM1}, \cite{DZ1}). In
particular, we show below (Proposition \ref{Prop_F_Mfg_F_Lh}) that every
such $g$-twisted $n$-point function is related to an $n$-point function for
the original VOSA but with a "shifted" Virasoro vector \cite{MN2}, \cite{DM}%
. This generalizes a similar result for partition functions found in \cite%
{DM} and allows us to apply Theorem \ref{Theorem_npt_rec} in order to
compute all such $g$-twisted $n$-point functions. The general relationship
at the operator level between these shifted and twisted formalisms is
discussed elsewhere \cite{TZ}.

A Heisenberg bosonic vector is an element $h\in V_{\bar{0},1}$ such that 
\cite{DM}

\begin{enumerate}
\item $h(0)$ is semisimple with real eigenvalues.

\item $h$ is a primary vector so that $L(n)h=0$ for all $n\geq 1$.

\item $h(n)h=0$ for all $n\geq 0$ except $n=1$ for which $h(1)h=\xi _{h}%
\mathbf{1}$ for some $\xi _{h}\in \mathbb{C}$.

\item $[h(m),h(n)]=\xi _{h}m\delta _{m,-n}$.
\end{enumerate}

\begin{remark}
\label{Remark_heisen}If the VOSA grading is non-negative and $V_{0}=\mathbb{C%
}\mathbf{1}$ then (2)-(4) follow automatically for all $h\in V_{\bar{0},1}$
from (\ref{aYcomm}).
\end{remark}

Given a Heisenberg vector $h$ then $h(0)$ generates a VOSA automorphism 
\begin{equation}
g=\exp (2\pi ih(0)).  \label{g_Heis}
\end{equation}%
The order of $g$ is finite iff the eigenvalues of $h(0)$ are rational and
otherwise is infinite. We can define \cite{DLM2} and construct a $g$-twisted
module in all cases as follows. We define \cite{Li} 
\begin{equation}
\Delta (h,z)=z^{h(0)}\exp \left( -\sum\limits_{n\geq 1}\frac{h(n)}{n}%
(-z)^{-n}\right) .  \label{Delta}
\end{equation}%
Noting $\Delta (h,z)^{-1}=\Delta (-h,z)$ one finds 
\begin{equation}
\Delta (h,z)Y(v,z_{0})\Delta (-h,z)=Y(\Delta (h,z+z_{0})v,z_{0}).
\label{delpro}
\end{equation}%
This leads to:

\begin{proposition}[\protect\cite{Li}]
\label{Proposition_YMh}Let $(M,Y_{M})$ be a $V$-module. Defining 
\begin{equation}
Y_{g}(v,z)=Y_{M}(\Delta (-h,z)v,z),  \label{Y_Mg}
\end{equation}%
for all $v\in V$ then $(M,Y_{g})$ is a $g$-twisted $V$-module\footnote{%
Note that we apply the definition of $g$-twisted module of ref. \cite{Li}
which corresponds to a $g^{-1}$-twisted module in refs. \cite{DLM1} and \cite%
{DM}}.
\end{proposition}

Note for $v=\omega $ we find $\Delta (-h,z)\omega =\omega -hz^{-1}+\xi
_{h}z^{-2}/2$ so that the $(M,Y_{g})$ grading is determined by 
\begin{equation}
L_{g}(0)=L(0)-h(0)+\frac{\xi _{h}}{2}.  \label{L0Mg}
\end{equation}%
We define the orbifold $g$-twisted $n$-point function for any automorphism $%
f $ commuting with $g$ and $\sigma $ by%
\begin{equation}
F_{M}((f,g);v_{1},\ldots ,v_{n};\tau )=\mathrm{STr}_{M}\left(
f\;Y_{g}(q_{1}^{L(0)}v_{1},q_{1})\ldots
Y_{g}(q_{n}^{L(0)}v_{n},q_{n})q^{L_{g}(0)-c/24}\right) .
\label{npointfunction_M_g}
\end{equation}%
We denote the orbifold $g$-twisted partition function by $Z_{M}((f,g),\tau )$%
.

For each Heisenberg element $h$ we may also construct a VOSA $(V,Y,\mathbf{1}%
,\omega _{h})$ with the original vector space and vertex operators but using
a "shifted" conformal vector (\cite{MN2}, \cite{DM})%
\begin{equation}
\omega _{h}=\omega +h(-2)\mathbf{1.}  \label{omega_h}
\end{equation}%
With $Y(\omega _{h},z)=\sum_{n\in \mathbb{Z}}L_{h}(n)z^{-n-2}$ we find 
\begin{equation}
L_{h}(n)=L(n)-(n+1)h(n),  \label{L_h(n)}
\end{equation}%
and central charge%
\begin{equation}
c_{h}=c-12\xi _{h}.  \label{C_h}
\end{equation}%
In particular, $L_{h}(-1)=L(-1)$ and the grading is determined by 
\begin{equation}
L_{h}(0)=L(0)-h(0).  \label{Lh(0)}
\end{equation}

We denote the partition function for a $V$-module $M$ with a $h$-shifted $%
L_{h}(0)$ by $Z_{M,h}(\tau )$. Following (\ref{npointfunction}) the shifted $%
n$-point function is denoted by 
\begin{equation}
F_{M,h}(f;v_{1},\ldots ,v_{n};\tau )=\mathrm{STr}_{M}\left(
fY(q_{1}^{L_{h}(0)}v_{1},q_{1})\ldots
Y(q_{n}^{L_{h}(0)}v_{n},q_{n})q^{L_{h}(0)-c_{h}/24}\right) ,
\label{npointfunction_h}
\end{equation}%
where $f$ commutes with $g$ and $\sigma $. We denote the $h$-shifted
partition function by $Z_{M,h}(f,\tau )$. Comparing (\ref{L0Mg}) and (\ref%
{Lh(0)}) we see that 
\begin{equation*}
L_{g}(0)-\frac{c}{24}=L_{h}(0)-\frac{c_{h}}{24},
\end{equation*}%
so that (\cite{DM}) 
\begin{equation}
Z_{M}((1,g),\tau )=Z_{M,h}(1,\tau ).  \label{twpart}
\end{equation}

\medskip This relationship can be generalized to relate all orbifold $g$%
-twisted $n$-point functions to $h$-shifted $n$-point functions as follows:

\begin{proposition}
\label{Prop_F_Mfg_F_Lh} Let $M$ be a module for $V$ and let $g=\exp (2\pi
ih(0))$ be generated by a Heisenberg state $h$. Then the $n$-point function
for the orbifold $g$-twisted and the untwisted $n$-point function for $M$
with shifted $L_{h}(0)$-vertex operators are related as follows:%
\begin{equation}
F_{M}((f,g);v_{1},\ldots ,v_{n};\tau )=F_{M,h}(f;Uv_{1},\ldots ,Uv_{n};\tau
),  \label{FMFh}
\end{equation}%
where $U=\Delta (-h,1)=\exp \left( \sum\limits_{n\geq 1}\frac{h(n)}{n}%
(-1)^{n}\right) $ and $f$ commutes with $g$ and $\sigma $.
\end{proposition}

\noindent \textbf{Proof. } First we prove 
\begin{equation}
\Delta (-h,q_{z})\;q_{z}^{L(0)}=q_{z}^{L_{h}(0)}\;U.  \label{urelation}
\end{equation}%
From (\ref{Delta}) one finds using $[L(0),h(n)]=-nh(n)$ that 
\begin{eqnarray*}
\Delta (-h,q_{z})q_{z}^{L(0)} &=&q_{z}^{-h(0)}\exp \left( \sum\limits_{n>0}%
\frac{h(n)}{n}(-q_{z})^{-n}\right) q_{z}^{L(0)} \\
&=&q_{z}^{-h(0)}q_{z}^{L(0)}\exp \left( \exp (\mathrm{ad}_{-zL(0)})\sum%
\limits_{n>0}\frac{h(n)}{n}(-q_{z})^{-n}\right) \\
&=&q_{z}^{L_{h}(0)}\exp \left( \sum\limits_{n>0}\frac{h(n)}{n}%
(-q_{z})^{-n}q_{z}^{n}\right) =q_{z}^{L_{h}(0)}U.
\end{eqnarray*}%
Therefore from (\ref{Y_Mg}) 
\begin{equation*}
Y_{g}(q_{z}^{L(0)}v,q_{z})=Y_{M}(\Delta
(-h,q_{z})q_{z}^{L(0)}v,q_{z})=Y_{M}(q_{z}^{L_{h}(0)}Uv,q_{z}),
\end{equation*}%
Thus the LHS of (\ref{FMFh}) is 
\begin{eqnarray*}
&&\mathrm{STr}_{M}\left( fY_{g}(q_{1}^{L(0)}v_{1},q_{1})\ldots
Y_{g}(q_{n}^{L(0)}v_{n},q_{n})q^{L_{g}(0)-c/24}\right) \\
&=&\mathrm{STr}_{M}\left( fY_{M}(q_{1}^{L_{h}(0)}Uv_{1},q_{1})\ldots
Y_{M}(q_{n}^{L_{h}(0)}Uv_{n},q_{n})q^{L_{h}(0)-c_{h}/24}\right) \\
&=&F_{M,h}(f;Uv_{1},\ldots ,Uv_{n};\tau ).
\end{eqnarray*}
\hfill $\square$

\medskip We conclude this section by showing that $U$ maps between the
square bracket vertex operators (\ref{Ysquare}) of the original and shifted
VOSAs. We let 
\begin{equation}
Y[v,z]_{h}=Y(q_{z}^{L_{h}(0)}v,q_{z}-1),  \label{squareh}
\end{equation}%
denote a square bracket vertex operator in the $h$-shifted VOSA. We then have

\begin{lemma}
\label{Lemma_Square}For $v\in V$ we have 
\begin{equation}
UY[v,z]U^{-1}=Y[Uv,z]_{h}.  \label{sqbraprop}
\end{equation}
\end{lemma}

\noindent \textbf{Proof. } Using associativity, (\ref{Y_Mg}) and (\ref%
{urelation}) we obtain 
\begin{eqnarray}
&&Y_{g}(q_{1}^{L(0)}v_{1},q_{1})Y_{g}(q_{2}^{L(0)}v_{2},q_{2})  \notag \\
&=&Y_{g}(Y(q_{1}^{L(0)}v_{1},q_{1}-q_{2})\;q_{2}^{L(0)}v_{2},q_{2})  \notag
\\
&=&Y_{g}(q_{2}^{L(0)}Y[v_{1},z_{12}]v_{2},q_{2})  \notag \\
&=&Y(q_{2}^{L_{h}(0)}UY[v_{1},z_{12}]v_{2},q_{2}).  \notag
\end{eqnarray}

On the other hand 
\begin{eqnarray*}
&&Y_{g}(q_{1}^{L(0)}v_{1},q_{1})Y_{g}(q_{2}^{L(0)}v_{2},q_{2}) \\
&=&Y(q_{1}^{L_{h}(0)}Uv_{1},q_{1})Y(q_{2}^{L_{h}(0)}Uv_{2},q_{2}) \\
&=&Y(q_{2}^{L_{h}(0)}Y[Uv_{1},z_{12}]_{h}\;Uv_{2},q_{2}).
\end{eqnarray*}%
Hence the result follows. \hfill $\square$

\section{\label{Section_Rank_One}Rank One Fermion VOSA}

We begin with the example of the rank one "Neveu-Schwarz sector" fermion
VOSA $V=V(H,\mathbb{Z}+\frac{1}{2})$ generated by one fermion \cite{FFR}, 
\cite{Li}. This is a $\frac{1}{2}\mathbb{Z}$ graded VOSA with $H=\mathbb{C}%
\psi $ for a fermion vector $\psi $ of parity $1$ and modes obeying%
\begin{equation}
\lbrack \psi (m),\psi (n)]=\psi (m)\psi (n)+\psi (n)\psi (m)=\delta
_{m+n+1,0}.  \label{Ferm_Com}
\end{equation}%
The superspace $V$ is spanned by Fock vectors of the form 
\begin{equation}
\psi (-k_{1})\psi (-k_{2})\ldots \psi (-k_{m})\mathbf{1},
\label{Fermion_Fock}
\end{equation}%
for integers $1\leq k_{1}<k_{2}<\ldots k_{m}$ with $\psi (k)\mathbf{1}=0$
for all $k\geq 0$ so that $V$ is generated by $Y(\psi ,z)$. The conformal
vector is $\omega =$ $\frac{1}{2}\psi (-2)\psi (-1)\mathbf{1}$ of central
charge $c=\frac{1}{2}$ for which the Fock vector (\ref{Fermion_Fock}) has $%
L(0)$ weight $\sum_{1\leq i\leq m}(k_{i}-\frac{1}{2})\in \frac{1}{2}\mathbb{Z%
}$. In particular, $wt(\psi )=\frac{1}{2}$. The partition function is 
\begin{equation}
Z_{V}(\tau )=\mathrm{STr}_{V}(q^{L(0)-\frac{1}{48}})=q^{-\frac{1}{48}%
}\prod_{n\geq 0}(1-q^{n+\frac{1}{2}})=\frac{\eta (\frac{1}{2}\tau )}{\eta
(\tau )},  \label{Zq_1Fermion}
\end{equation}%
whereas for $g=\sigma $ of (\ref{sigma}) we find%
\begin{equation}
Z_{V}(\sigma ,\tau )=\mathrm{STr}_{V}(\sigma q^{L(0)-\frac{1}{48}})=q^{-%
\frac{1}{48}}\prod_{n\geq 0}(1+q^{n+\frac{1}{2}})=\frac{\eta (\tau )^{2}}{%
\eta (2\tau )\eta (\frac{1}{2}\tau )}.  \label{Zq_1Ferm_sig}
\end{equation}

Let us next introduce the $n$-point function (\ref{npointfunction}) for $V$
where $v_{i}=\psi $ for all $i=1,\ldots n$: 
\begin{equation}
G_{n}(g;z_{1},\ldots ,z_{n};\tau )=F_{V}(g;(\psi ,z_{1}),\ldots ,(\psi
,z_{n});\tau ),  \label{genfu}
\end{equation}%
which we will refer to as the \emph{generating function}. We use the
recursion formula (\ref{nptrec}) of Theorem \ref{Theorem_npt_rec} to compute 
$G_{n}$. Since $wt(\psi )=\frac{1}{2}$ we have $\phi =-1$ from (\ref{phi})
and $\theta =1$ for $g=1$ and $\theta =-1$ for $g=\sigma $ from (\ref{theta}%
). For $n=1$, $G_{1}(g;z_{1};\tau )=Z_{V}(g,\psi ,\tau )=0$ since $o(\psi
)=0 $. For $n=2$, (\ref{nptrec}) implies 
\begin{equation*}
G_{2}(g;z_{1},z_{2};\tau )=0+\sum\limits_{m\geq 0}P_{m+1}\left[ 
\begin{array}{c}
\theta \\ 
-1%
\end{array}%
\right] (z_{12},\tau )F_{V}(g;\psi \lbrack m]\psi ;\tau ).
\end{equation*}%
Passing to the square bracket formalism (\ref{Ysquare}) we find the same
fermion commutator algebra as (\ref{Ferm_Com}) obtains, namely 
\begin{equation}
\lbrack \psi \lbrack m],\psi \lbrack n]]=\delta _{m+n+1,0}.
\label{Ferm_Com_Square}
\end{equation}%
Thus it follows that $\psi \lbrack m]\psi =\delta _{m,0}\mathbf{1}$ giving 
\begin{equation}
G_{2}(g;z_{1},z_{2};\tau )=P_{1}\left[ 
\begin{array}{c}
\theta \\ 
-1%
\end{array}%
\right] (z_{12},\tau )Z_{V}(g,\tau ).  \label{G2fermion}
\end{equation}

We may similarly compute $G_{n}$ for all $n$ by repeated application of (\ref%
{nptrec}). It is easy to see that $G_{n}=0$ for $n$ odd. For $n$ even $G_{n}$
is expressed in terms of a Pfaffian which is totally antisymmetric in $z_{i}$
as expected from Lemma \ref{Lemma_Perm_Period} (ii). Let us first recall the
definition of the Pfaffian of an anti-symmetric matrix $\mathbf{M}=(M(i,j))$
of even dimension $2m$ given by 
\begin{equation}
\mathrm{Pf}(\mathbf{M})=\sum_{\Pi }\varepsilon _{i_{1}j_{1}\ldots
i_{m}j_{m}}M(i_{1},j_{1})M(i_{2},j_{2})\ldots M(i_{m},j_{m}),  \label{Pfaff}
\end{equation}%
where the sum is taken over the set of all partitions $\Pi $ of $%
\{1,2,\ldots ,2m\}$ into pairs with elements%
\begin{equation*}
\{(i_{1},j_{1}),(i_{2},j_{2})\ldots (i_{m},j_{m})\},
\end{equation*}%
for $i_{k}<j_{k}$ and $i_{1}<i_{2}<\ldots i_{m}$ and where $\varepsilon
_{i_{1}j_{1}\ldots i_{m}j_{m}}$ is the Levi-Civita symbol. We also note that 
\begin{equation*}
\mathrm{Pf}(\mathbf{M})=\sqrt{\det \mathbf{M}}.
\end{equation*}%
We then obtain:

\begin{proposition}
\label{PropGnPfaff} For $n$ even and $g=1$ or $\sigma $ we have 
\begin{equation}
G_{n}(g;z_{1},\ldots ,z_{n};\tau )=\mathrm{Pf}(\mathbf{P})Z_{V}(g,\tau ),
\label{GnPfaff}
\end{equation}%
where $\mathbf{P}$ denotes the anti-symmetric $n\times n$ matrix with
components 
\begin{equation}
\mathbf{P}(i,j)=P_{1}\left[ 
\begin{array}{c}
\theta \\ 
-1%
\end{array}%
\right] (z_{ij},\tau ),\quad (1\leq i\neq j\leq n),  \label{Ptheta}
\end{equation}%
for $z_{ij}=z_{i}-z_{j}$ with $\theta =1$ for $g=1$ and $\theta =-1$ for $%
g=\sigma $.
\end{proposition}

\noindent \textbf{Proof. } We first note that $\mathbf{P}$ is anti-symmetric
from (\ref{Pkodd}) since $\theta =\pm 1$. We prove the result by induction.
For $n=2$ the result is given in (\ref{G2fermion}). For general $n$ we apply
(\ref{nptrec}) to obtain 
\begin{eqnarray*}
G_{n}(g;z_{1},\ldots z_{n};\tau ) &=&\sum\limits_{r=2}^{n}(-1)^{r}P_{1}\left[
\begin{array}{c}
\theta \\ 
-1%
\end{array}%
\right] (z_{1r},\tau )G_{n-2}(g;z_{2},\ldots ,\hat{z}_{r},\ldots z_{n};\tau )
\\
&=&\sum\limits_{r=2}^{n}(-1)^{r}\mathbf{P}(1,r)\mathrm{Pf}(\mathbf{\hat{P}}%
)Z_{V}(g,\tau ),
\end{eqnarray*}%
where $\hat{z}_{r}$ is deleted and $\mathbf{\hat{P}}$ is the "cofactor"
matrix obtained by deleting the $1^{st}$ and $r^{th}$ rows and columns of $%
\mathbf{P}$. The result (\ref{Ptheta}) follows from the definition (\ref%
{Pfaff}).$\square $

\medskip $G_{n}$ enjoys the following analytic properties following Remark %
\ref{Remark_Pk} (ii):

\begin{corollary}
\label{Corollary_GnPfaff}$G_{n}$ is an analytic function in $z_{i}$ and
converges absolutely and uniformly on compact subsets of the domain $%
\left\vert q\right\vert <\left\vert q_{z_{ij}}\right\vert <1$ for all $%
z_{ij}=z_{i}-z_{j}$ with $i\neq j$. $\square $
\end{corollary}

We now show that all $n$-point functions can be computed from $G_{n}$.
Consider a $V$ basis of square bracket Fock vectors denoted by 
\begin{equation}
\Psi \lbrack -\mathbf{k}]=\psi \lbrack -k_{1}]\psi \lbrack -k_{2}]\ldots
\psi \lbrack -k_{m}]\mathbf{1},  \label{Fermion_Fock_Square}
\end{equation}%
where $\mathbf{k}=k_{1},k_{2},\ldots ,k_{m}$ for integers $1\leq
k_{1}<k_{2}<\ldots k_{m}$. We will determine an explicit formula for all $n$%
-point functions for such Fock vectors. Thus the 1-point function $%
Z_{V}(g,\Psi \lbrack -\mathbf{k}],\tau )$ is the coefficient of $%
\prod_{i=1}^{m}z_{i}^{k_{i}-1}$ in $G_{n}$ since 
\begin{eqnarray*}
G_{n}(g;z_{1},\ldots ,z_{m};\tau ) &=&Z_{V}(g,Y[\psi ,z_{1}]\ldots Y[\psi
,z_{m}]\mathbf{1},\tau ) \\
&=&\sum_{k_{1},\ldots k_{m}\in \mathbb{Z}}Z_{V}(g,\psi \lbrack -k_{1}]\ldots
\psi \lbrack -k_{m}]\mathbf{1},\tau )z_{1}^{k_{1}-1}\ldots z_{n}^{k_{m}-1}.
\end{eqnarray*}%
Examining (\ref{GnPfaff}) we can explicitly find this coefficient from the
expansion of $P_{1}\left[ {%
\begin{array}{c}
\theta \\ 
-1%
\end{array}%
}\right] (z_{ij},\tau )$ given in (\ref{P1_C_expansion}). It follows that $%
Z_{V}(g,\Psi \lbrack -\mathbf{k}],\tau )=0$ for $m$ odd whereas for $m$ even 
\begin{equation}
Z_{V}(g,\Psi \lbrack -\mathbf{k}],\tau )=\mathrm{Pf}(\mathbf{C})Z_{V}(g,\tau
),  \label{Z1ptFockC}
\end{equation}%
where $\mathbf{C}$ denotes the antisymmetric $m\times m$ matrix with $(i,j)$%
-entry%
\begin{equation*}
\mathbf{C}(i,j)=C\left[ 
\begin{array}{c}
\theta \\ 
-1%
\end{array}%
\right] (k_{i},k_{j},\tau ),
\end{equation*}%
(cf. (\ref{Ckldef})). $\mathbf{C}$ is antisymmetric from (\ref{Cklodd})
since $\theta =\pm 1$.

We may similarly derive an expression for an arbitrary two-point function $%
F_{V}((\Psi \lbrack -\mathbf{k}^{(1)}],z_{1}),(\Psi \lbrack -\mathbf{k}%
^{(2)}],z_{2});g;\tau )$ for $\mathbf{k}^{(1)}=k_{1}^{(1)},\ldots
k_{m_{1}}^{(1)}$ and $\mathbf{k}^{(2)}=k_{1}^{(2)},\ldots k_{m_{2}}^{(2)}$.
First consider the one-point function 
\begin{equation}
Z_{V}(g,Y[Y[\psi ,x_{1}]\ldots Y[\psi ,x_{m_{1}}]\mathbf{1},z_{1}].Y[Y[\psi
,y_{1}]\ldots Y[\psi ,y_{m_{2}}]\mathbf{1},z_{2}]\mathbf{1},\tau ).
\label{FNtwoptgen}
\end{equation}%
$F_{V}(g;(\Psi \lbrack -\mathbf{k}^{(1)}],z_{1}),(\Psi \lbrack -\mathbf{k}%
^{(2)}],z_{2});\tau )$ is the coefficient of $\prod_{i=1}^{m_{1}}%
\prod_{j=1}^{m_{2}}x_{i}^{k_{i}^{(1)}-1}y_{j}^{k_{j}^{(2)}-1}$ in (\ref%
{FNtwoptgen}). By associativity (e.g. \cite{FHL}) and using $Y[\mathbf{1},z]=%
\mathrm{Id}_{V}$ we find (\ref{FNtwoptgen}) can be expressed as 
\begin{eqnarray*}
&&Z_{V}(g,Y[\psi ,x_{1}+z_{1}]\ldots Y[\psi ,x_{m_{1}}+z_{1}].Y[\psi
,y_{1}+z_{2}]\ldots Y[\psi ,y_{m_{2}}+z_{2}]\mathbf{1},\tau ) \\
&=&G_{n}(g;x_{1}+z_{1},\ldots ,x_{m_{1}}+z_{1},y_{1}+z_{2},\ldots
,y_{m_{2}}+z_{2};\tau ).
\end{eqnarray*}%
The coefficient of $\prod_{i=1}^{m_{1}}%
\prod_{j=1}^{m_{2}}x_{i}^{k_{i}^{(1)}-1}y_{j}^{k_{j}^{(2)}-1}$ can then be
extracted from the expansions (\ref{P1_C_expansion}) and (\ref%
{P1_D_expansion}). Thus the two point function vanishes for $m_{1}+m_{2}$
odd, whereas for $m_{1}+m_{2}$ even 
\begin{equation}
F_{V}(g;(\Psi \lbrack -\mathbf{k}^{(1)}],z_{1}),(\Psi \lbrack -\mathbf{k}%
^{(2)}],z_{2});\tau )=\mathrm{Pf}(\mathbf{M})Z(g,\tau ),  \label{F2ptFockCD}
\end{equation}%
where $\mathbf{M}$ is the antisymmetric $(m_{1}+m_{2})\times (m_{1}+m_{2})$
block matrix 
\begin{equation*}
\mathbf{M=}\left( 
\begin{array}{cc}
\mathbf{C}^{(11)} & \mathbf{D}^{(12)} \\ 
\mathbf{D}^{(21)} & \mathbf{C}^{(22)}%
\end{array}%
\right) ,
\end{equation*}%
where for $a,b\in \{1,2\}$ 
\begin{eqnarray}
\mathbf{C}^{(aa)}(i,j) &=&C\left[ 
\begin{array}{c}
\theta \\ 
-1%
\end{array}%
\right] (k_{i}^{(a)},k_{j}^{(a)},\tau ),\quad (1\leq i,j\leq m_{a}),  \notag
\\
\mathbf{D}^{(ab)}(i,j) &=&D\left[ 
\begin{array}{c}
\theta \\ 
-1%
\end{array}%
\right] (k_{i}^{(a)},k_{j}^{(b)},\tau ,z_{a}-z_{b}),\quad (1\leq i\leq
m_{a},1\leq j\leq m_{b}),  \notag \\
&&  \label{CDmatrices}
\end{eqnarray}%
(using (\ref{Dkldef})). $\mathbf{M}$ is antisymmetric from (\ref{Cklodd})
and (\ref{Dklodd}).

\medskip In a similar fashion we are lead to the general result:

\begin{proposition}
\label{Prop_Ferm_1_npt} Let $\Psi \lbrack -\mathbf{k}^{(a)}]$ for $a=1\ldots
n$ be $n$ Fock vectors for $\mathbf{k}^{(a)}=k_{1}^{(a)},\ldots
k_{m_{a}}^{(a)}$. Then the $n$-point function vanishes for odd $%
\sum_{a}m_{a} $ and for $\sum_{a}m_{a}$ even is given by 
\begin{equation}
F_{V}(g;(\Psi \lbrack -\mathbf{k}^{(1)}],z_{1}),\ldots (\Psi \lbrack -%
\mathbf{k}^{(n)}],z_{n});\tau )=\mathrm{Pf}(\mathbf{M})Z(g,\tau ),
\label{FnptFock}
\end{equation}%
where $\mathbf{M}$ is the antisymmetric block matrix 
\begin{equation*}
\mathbf{M=}\left( 
\begin{array}{cccc}
\mathbf{C}^{(11)} & \mathbf{D}^{(12)} & \ldots & \mathbf{D}^{(1n)} \\ 
\mathbf{D}^{(21)} & \mathbf{C}^{(22)} &  &  \\ 
\vdots &  & \ddots &  \\ 
\mathbf{D}^{(n1)} & \ldots &  & \mathbf{C}^{(nn)}%
\end{array}%
\right) ,
\end{equation*}%
with $\mathbf{C}^{(aa)}$ and $\mathbf{D}^{(ab)}$ of (\ref{CDmatrices}). (\ref%
{FnptFock}) is an analytic function in $z_{i}$ and converges absolutely and
uniformly on compact subsets of the domain $\left\vert q\right\vert
<\left\vert q_{z_{ij}}\right\vert <1$ for all $z_{ij}=z_{i}-z_{j}$ with $%
i\neq j$. $\square $
\end{proposition}

We have also established

\begin{proposition}
\label{Gen_rankone} $G_{n}(g;z_{1},\ldots ,z_{n};\tau )$ is a generating
function for all $n$-point functions. $\square $
\end{proposition}

We conclude this section by noting that we may also consider the "Ramond
sector" $\sigma $-twisted module $V(H,\mathbb{Z)}$ for $V(H,\mathbb{Z}+\frac{%
1}{2})$. This is discussed in detail in \cite{FFR}, \cite{Li}, \cite{DZ1}, 
\cite{DZ2}. $V(H,\mathbb{Z)}$ decomposes into two irreducible $\sigma $%
-twisted modules which are interchanged under the induced action of $\sigma $%
. For either irreducible $\sigma $-twisted module $M_{\sigma }$ the
partition function is 
\begin{eqnarray*}
\mathrm{STr}_{M_{\sigma }}(q^{L(0)-\frac{1}{48}}) &=&0, \\
\mathrm{STr}_{M_{\sigma }}(\sigma q^{L(0)-\frac{1}{48}}) &=&q^{\frac{1}{48}%
}\prod_{n\geq 0}(1+q^{n})=\frac{\eta (2\tau )}{\eta (\tau )}.
\end{eqnarray*}%
We may similarly consider the generator of all $\sigma $-twisted $n$-point
functions defined by 
\begin{equation*}
G_{M_{\sigma },n}(g;z_{1},\ldots ,z_{n};\tau )=F_{M_{\sigma }}(g;(\psi
,z_{1}),\ldots ,(\psi ,z_{n});\tau ),
\end{equation*}%
for $g=1$ or $\sigma $. This vanishes for all $n$ for $g=1$ and for $n$ odd
for $g=\sigma $. By applying a VOSA orbifold Zhu reduction formula of ref. 
\cite{DZ1} we find as in Proposition \ref{PropGnPfaff} that

\begin{proposition}
\label{PropGnPfaffsig} For $n$ even we have 
\begin{equation}
G_{M_{\sigma },n}(\sigma ;z_{1},\ldots ,z_{n};\tau )=\mathrm{Pf}(P_{1}\left[ 
\begin{array}{c}
-1 \\ 
1%
\end{array}%
\right] (z_{ij},\tau ))\frac{\eta (2\tau )}{\eta (\tau )},  \label{Gn_Msigma}
\end{equation}%
for $z_{ij}=z_{i}-z_{j}$. $\square $
\end{proposition}

One can similarly describe analytic properties as in Corollary \ref%
{Corollary_GnPfaff} and determine all $\sigma $-twisted $n$-point functions
by expanding this generating function along the same lines as Proposition %
\ref{Prop_Ferm_1_npt}, though we do not carry this out here.

\section{\label{Section_Rank_Two}Rank Two Fermion VOSA}

\medskip

\subsection{$h$-Shifted and orbifold $g$-Twisted $n$-Point Functions}

In this section we consider the rank two fermion VOSA formed from the tensor
product of two copies of the rank one fermion VOSA and hence is generated by
two free fermions $\psi _{1}=\psi \otimes \mathbf{1}$ and $\psi _{2}=\mathbf{%
1}\otimes \psi $. We may therefore compute all the untwisted and $\sigma $%
-twisted $n$-point functions based on the last section. However, as is well
known, this VOSA contains a bosonic Heisenberg state $h=\alpha \psi \otimes
\psi $ (for $\alpha \in \mathbb{C}$) and we will compute all $h$-shifted and 
$g$-twisted $n$-point functions where $g$ is generated by $h$ as discussed
in Section \ref{Section_Heisenberg}.

It is convenient to introduce the off-diagonal basis $\psi ^{\pm }=\frac{1}{%
\sqrt{2}}\left( \psi _{1}\pm i\psi _{2}\right) $, where $\psi ^{\pm }$-modes
obey the commutation relations 
\begin{equation}
\lbrack \psi ^{+}(m),\psi ^{-}(n)]=\delta _{m,-n-1},\quad \lbrack \psi ^{\pm
}(m),\psi ^{\pm }(n)]=0,  \label{psiplus_minus_comm}
\end{equation}
The VOSA $V$ is generated by $Y(\psi ^{\pm },z)=\sum\limits_{n\in \mathbb{Z}%
}\psi ^{\pm }(n)z^{-n-1}$ where the vector space $V$ is a Fock space with
basis vectors of the form 
\begin{equation}
\psi ^{+}(-k_{1})...\psi ^{+}(-k_{s})\psi ^{-}(-l_{1})...\psi ^{-}(-l_{t})%
\mathbf{1},  \label{Fockstate}
\end{equation}%
for $1\leq k_{1}<k_{2}<\ldots k_{s}$ and $1\leq l_{1}<l_{2}<\ldots l_{t}$
with $\psi ^{\pm }(k)\mathbf{1}=0$ for all $k\geq 0$. We define the
conformal vector to be 
\begin{equation}
\omega =\frac{1}{2}[\psi ^{+}(-2)\psi ^{-}(-1)\mathbf{+}\psi ^{-}(-2)\psi
^{+}(-1)]\mathbf{1},  \label{omega_psi}
\end{equation}%
whose modes generate a Virasoro algebra of central charge 1. Then $\psi
^{\pm }$ has $L(0)$-weight $\frac{1}{2}$ and the Fock state (\ref{Fockstate}%
) has weight $\sum_{1\leq i\leq s}(k_{i}-\frac{1}{2})+\sum_{1\leq j\leq
t}(l_{j}-\frac{1}{2})$.

\medskip The weight $1$ parity zero space is $V_{\bar{0},1}=\mathbb{C}a$ for
(normalized) Heisenberg bosonic vector 
\begin{equation}
a=\psi ^{+}(-1)\psi ^{-}(-1)\mathbf{1,}  \label{adef}
\end{equation}%
with modes obeying 
\begin{equation*}
\lbrack a(m),a(n)]=m\delta _{m,-n},
\end{equation*}%
and $\omega $ of (\ref{omega_psi}) is nothing but the standard Heisenberg
VOA conformal vector 
\begin{equation*}
\omega =\frac{1}{2}a(-1)^{2}\mathbf{1}.
\end{equation*}

Following Section \ref{Section_Heisenberg} we define a one parameter family
of Heisenberg vectors 
\begin{equation}
h=\kappa a,\quad \kappa \in \mathbb{R},  \label{h_kappa}
\end{equation}%
for which $\xi _{h}=\kappa ^{2}$. The shifted conformal vector (\ref{omega_h}%
) is then $\omega _{h}=\omega +\kappa a(-2)\mathbf{1}$ with central charge $%
c_{h}=1-12\kappa ^{2}$ from (\ref{C_h}). Then $\psi ^{\pm }$ has $%
L_{h}(0)=L(0)-\kappa a(0)$ weight $wt_{h}(\psi ^{\pm })=$ $\frac{1}{2}\mp
\kappa $ and the Fock state (\ref{Fockstate}) has $L_{h}(0)$ weight $%
\sum_{1\leq i\leq s}(k_{i}-\frac{1}{2}-\kappa )+\sum_{1\leq j\leq t}(l_{j}-%
\frac{1}{2}+\kappa )$.

\medskip Noting that $\sigma =e^{\pi ia(0)}$ and following Section \ref%
{Section_Heisenberg}, we can construct a $\sigma g$-twisted module for $%
\sigma g=e^{2\pi ih(0)}$ so that 
\begin{equation}
g=e^{2\pi i\beta a(0)},  \label{g_beta}
\end{equation}%
for real $\beta $ where 
\begin{equation}
\beta =\kappa -\frac{1}{2}.  \label{beta_kappa}
\end{equation}%
We also define $\phi \in U(1)\ $by%
\begin{equation}
\phi =\exp (2\pi iwt_{h}(\psi ^{+}))=e^{-2\pi i\beta }.  \label{phi_beta}
\end{equation}%
Introduce the automorphism 
\begin{equation}
f=e^{2\pi i\alpha a(0)},\quad \alpha \in \mathbb{R},  \label{f_alpha}
\end{equation}%
which commutes with $g,\sigma $. Then $f\psi ^{\pm }=\theta ^{\mp 1}\psi
^{\pm }$ for%
\begin{equation}
\theta =e^{-2\pi i\alpha }\in U(1).  \label{theta_alpha}
\end{equation}%
Finally, we denote the orbifold $\sigma g$-twisted trace by 
\begin{equation*}
Z_{V}\left[ 
\begin{array}{c}
f \\ 
g%
\end{array}%
\right] (\tau )=Z_{V}((f,\sigma g),\tau ).
\end{equation*}%
We find using Proposition \ref{Prop_F_Mfg_F_Lh} that 
\begin{equation}
Z_{V,h}(f,\tau )=Z_{V}\left[ 
\begin{array}{c}
f \\ 
g%
\end{array}%
\right] (\tau )=q^{\kappa ^{2}/2-1/24}\prod_{l\geq 1}(1-\theta ^{-1}q^{l-%
\frac{1}{2}-\kappa })(1-\theta q^{l-\frac{1}{2}+\kappa }).  \label{Zparth}
\end{equation}%
Note that $Z_{V,h}(f,\tau )=0$ for $(\theta ,\phi )=(1,1)$, i.e. $(\alpha
,\beta ) \equiv (0,0) \ (\mbox{mod} \ \mathbb{Z}$).

\begin{remark}
\label{Remark_part_alpha}The RHS of (\ref{Zparth}) is related to a theta
series via the Jacobi triple product formula as briefly reviewed below in
Section \ref{Boson}. Hence $Z_{V,h}(f,\tau )$ depends on $\alpha (\mathrm{{%
mod} \ \mathbb{Z})}$ and $\beta (\mathrm{{mod} \ \mathbb{Z})}$ up to an
overall $\alpha$-dependent constant.
\end{remark}

We next consider general $\sigma g$-twisted and $h$-shifted $n$-point
functions which are related via Proposition \ref{Prop_F_Mfg_F_Lh}. As in the
rank one case, it is sufficient to consider $n$-point functions for the
generating states $\psi ^{\pm }$ only. To this end we define the $h$-shifted
VOSA $n$-point generating function 
\begin{eqnarray}
&&G_{2n,h}(f;x_{1},...,x_{n};y_{1},...,y_{n};\tau )  \notag \\
&=&F_{V,h}(f;(\psi ^{+},x_{1}),(\psi ^{-},y_{1}),...,(\psi ^{+},x_{n}),(\psi
^{-},y_{n});\tau ),  \label{F_VLh_gen}
\end{eqnarray}

\begin{remark}
\label{Remark_check}Note the choice of an alternating ordering of the
operators with respect to the $\pm $ superscript here.
\end{remark}

We can also define a $\sigma g$-twisted $n$-point function denoted by%
\begin{equation*}
F_{V}\left[ 
\begin{array}{c}
f \\ 
g%
\end{array}%
\right] ((v_{1},z_{1})...,(v_{n},z_{n});\tau )=F_{V}((f,\sigma
g);(v_{1},z_{1})...,(v_{n},z_{n});\tau ),
\end{equation*}%
with generating function 
\begin{eqnarray*}
&&G_{2n}\left[ 
\begin{array}{c}
f \\ 
g%
\end{array}%
\right] (x_{1},...,x_{n};y_{1},...,y_{n};\tau ) \\
&=&F_{V}((f,\sigma g);(\psi ^{+},x_{1}),(\psi ^{-},y_{1}),...,(\psi
^{+},x_{n}),(\psi ^{-},y_{n});\tau ).
\end{eqnarray*}%
Then noting that $U$ $\psi ^{\pm }=\psi ^{\pm }$ and applying Proposition %
\ref{Prop_F_Mfg_F_Lh} we find

\begin{lemma}
\label{lemma_Gn_h_g}%
\begin{equation*}
G_{2n}\left[ 
\begin{array}{c}
f \\ 
\sigma g%
\end{array}%
\right] (x_{1},...,x_{n};y_{1},...,y_{n};\tau
)=G_{2n,h}(f;x_{1},...,x_{n};y_{1},...,y_{n};\tau ).\quad \square
\end{equation*}
\end{lemma}

These generating functions are totally antisymmetric in $x_{i},y_{j}$ as
expected from Lemma \ref{Lemma_Perm_Period} (ii) and can be expressed in
terms of a determinant computed by means of our recursion formula (\ref%
{nptrec}). Due to the leading term on the RHS of (\ref{nptrec}), we consider
the cases $(\theta ,\phi )\neq (1,1)$ and $(\theta ,\phi )=(1,1)$ separately.

\subsection{$n$-Point Functions for $(\protect\theta ,\protect\phi )\neq
(1,1)$.}

\begin{proposition}
\label{prop_genVh} For $(\theta ,\phi )\neq (1,1)$ we have 
\begin{equation}
G_{2n,h}(f;x_{1},...,x_{n};y_{1},...,y_{n};\tau )=\mathrm{\det }\mathbf{P}%
.\;Z_{V,h}(f;\tau ),  \label{FMonepointgen11}
\end{equation}%
where $\mathbf{P}$ is the $n\times n$ matrix: 
\begin{equation}
\mathbf{P}=\left( P_{1}\left[ 
\begin{array}{c}
\theta \\ 
\phi%
\end{array}%
\right] (x_{i}-y_{j},\tau )\right) ,\quad (1\leq i,j\leq n),  \label{Pmatrix}
\end{equation}%
with $\theta ,\phi $ of (\ref{theta_alpha}) and (\ref{phi_beta}).
Furthermore, $G_{2n,h}$ is an analytic function in $x_{i},y_{j}$ and
converges absolutely and uniformly on compact subsets of the domain $%
\left\vert q\right\vert <\left\vert q_{x_{i}-y_{j}}\right\vert <1$.
\end{proposition}

\noindent \textbf{Proof. }We apply Theorem \ref{Theorem_npt_rec} directly
with $e^{2\pi iwt(\psi ^{+})}=\phi $ and $f\psi ^{+}=e^{2\pi i\beta }\psi
^{+}=\theta ^{-1}\psi ^{+}$ of (\ref{phi_beta}) and (\ref{theta_alpha}). The
Zhu recursion formula (\ref{nptrec}) results in a determinant similarly to
the proof of Proposition \ref{PropGnPfaff}. The region of analyticity
follows as before. $\square $

\medskip In order to describe general $n$-point functions, first note that 
\begin{equation*}
\lbrack \psi ^{+}[m],\psi ^{-}[n]]=\delta _{m,-n-1},\quad \lbrack \psi ^{\pm
}[m],\psi ^{\pm }[n]]=0.
\end{equation*}%
Now introduce 
\begin{eqnarray}
\Psi &=&\Psi \lbrack -\mathbf{k};-\mathbf{l}]=\psi ^{+}[-k_{1}]...\psi
^{+}[-k_{s}]\psi ^{-}[-l_{1}]...\psi ^{-}[-l_{t}]\mathbf{1},
\label{Square_Fock} \\
\Psi _{h} &=&\Psi \lbrack -\mathbf{k};-\mathbf{l}]_{h}=\psi
^{+}[-k_{1}]_{h}...\psi ^{+}[-k_{s}]_{h}\psi ^{-}[-l_{1}]_{h}...\psi
^{-}[-l_{t}]_{h}\mathbf{1},  \label{Square_Fock_h}
\end{eqnarray}%
where $\mathbf{k}=k_{1},...,k_{s}$ and $\mathbf{l}=l_{1},...,l_{t}$; these
denote Fock vectors (\ref{Fockstate}) in the square bracket and $h$-shifted
square bracket formalisms respectively. From Lemma \ref{Lemma_Square} and
using $U$ $\psi ^{\pm }=\psi ^{\pm }$ we have%
\begin{equation*}
\Psi \lbrack -\mathbf{k};-\mathbf{l}]_{h}=U\Psi \lbrack -\mathbf{k};-\mathbf{%
l}].
\end{equation*}%
\ \ By expanding $G_{2n,h}$ appropriately and following the same approach
that lead to Proposition \ref{Prop_Ferm_1_npt}, we obtain a determinant
formula for every $n$-point function as follows:

\begin{proposition}
\label{prop_npt_Vh_g} Consider $n$ Fock vectors $\Psi ^{(a)}=\Psi ^{(a)}[-%
\mathbf{k}^{(a)};-\mathbf{l}^{(a)}]$ and $\Psi _{h}^{(a)}=\Psi ^{(a)}[-%
\mathbf{k}^{(a)};-\mathbf{l}^{(a)}]_{h}$ for $\mathbf{k}%
^{(a)}=k_{1}^{(a)},...k_{s_{a}}^{(a)}$ and $\mathbf{l}%
^{(a)}=l_{1}^{(a)},...l_{t_{a}}^{(a)}$ with $a=1\ldots n$. Then for $(\theta
,\phi )\neq (1,1)$ the corresponding $n$-point functions are non-vanishing
provided 
\begin{equation*}
\sum\limits_{a=1}^{n}\left( s_{a}-t_{a}\right) =0.
\end{equation*}%
In this case they are given by 
\begin{eqnarray}
&&F_{V}\left[ 
\begin{array}{c}
f \\ 
g%
\end{array}%
\right] ((\Psi ^{(1)},z_{1}),\ldots ,(\Psi ^{(n)},z_{n});\tau )  \notag \\
&=&F_{V,h}(f;(\Psi _{h}^{(1)},z_{1}),\ldots ,(\Psi _{h}^{(n)},z_{n});\tau
)=\epsilon \;\mathrm{\det }\mathbf{M.}\;Z_{V,h}(f;\tau ),  \label{Fnpt_rank2}
\end{eqnarray}%
where $\mathbf{M}$ is the block matrix%
\begin{equation*}
\mathbf{M}=\left( 
\begin{array}{ccc}
\mathbf{C}^{(11)} & \mathbf{D}^{(12)}\ldots & \mathbf{D}^{(1n)} \\ 
\mathbf{D}^{(21)} & \mathbf{C}^{(22)}\ldots & \mathbf{D}^{(2n)} \\ 
\vdots & \ddots & \vdots \\ 
\mathbf{D}^{(n1)} & \ldots & \mathbf{C}^{(nn)}%
\end{array}%
\right) ,
\end{equation*}%
with 
\begin{equation*}
\mathbf{C}^{(aa)}(i,j)=C\left[ 
\begin{array}{c}
\theta \\ 
\phi%
\end{array}%
\right] (k_{i}^{(a)},l_{j}^{(a)},\tau ),\quad (1\leq i\leq s_{a},1\leq j\leq
t_{a}),
\end{equation*}%
for $s_{a},t_{a}\geq 1$ with $1\leq a\leq n$ and 
\begin{equation*}
\mathbf{D}^{(ab)}(i,j)=D\left[ 
\begin{array}{c}
\theta \\ 
\phi%
\end{array}%
\right] (k_{i}^{(a)},l_{j}^{(b)},\tau ,z_{ab}),\quad (1\leq i\leq
s_{a},1\leq j\leq t_{b}),
\end{equation*}%
for $s_{a},t_{b}\geq 1$ with $1\leq a,b\leq n\ $and $a\neq b$. $\epsilon $
is the sign of the permutation associated with the reordering of $\psi ^{\pm
}$ to the alternating ordering of (\ref{F_VLh_gen}) following Remark \ref%
{Remark_check}. Furthermore, the $n$-point function (\ref{Fnpt_rank2}) is an
analytic function in $z_{a}$ and converges absolutely and uniformly on
compact subsets of the domain $\left\vert q\right\vert <\left\vert
q_{z_{ab}}\right\vert <1$. $\square $
\end{proposition}

\medskip \noindent \textbf{Example.} \label{Example}Consider the $n$-point
function for $n$ vectors $\Psi =a$ for $a=\psi ^{+}[-1]\psi ^{-}[-1]\mathbf{1%
}$ and $(\theta ,\phi )\neq (1,1)$. Then 
\begin{equation*}
F_{V}\left[ 
\begin{array}{c}
f \\ 
g%
\end{array}%
\right] ((a,z_{1}),\ldots ,(a,z_{n});\tau )=\det M.Z_{V}\left[ 
\begin{array}{c}
f \\ 
g%
\end{array}%
\right] (\tau ),
\end{equation*}%
for 
\begin{equation*}
\mathbf{M}=\left( 
\begin{array}{ccc}
-E_{1}\left[ {%
\begin{array}{c}
\theta \\ 
\phi%
\end{array}%
}\right] (\tau ) & P_{1}\left[ {%
\begin{array}{c}
\theta \\ 
\phi%
\end{array}%
}\right] (z_{12},\tau )\ldots & P_{1}\left[ {%
\begin{array}{c}
\theta \\ 
\phi%
\end{array}%
}\right] (z_{1n},\tau ) \\ 
P_{1}\left[ {%
\begin{array}{c}
\theta \\ 
\phi%
\end{array}%
}\right] (z_{21},\tau ) & -E_{1}\left[ {%
\begin{array}{c}
\theta \\ 
\phi%
\end{array}%
}\right] (\tau )\ldots & P_{1}\left[ {%
\begin{array}{c}
\theta \\ 
\phi%
\end{array}%
}\right] (z_{2n},\tau ) \\ 
\vdots & \ddots & \vdots \\ 
P_{1}\left[ {%
\begin{array}{c}
\theta \\ 
\phi%
\end{array}%
}\right] (z_{n1},\tau ) & \ldots & -E_{1}\left[ {%
\begin{array}{c}
\theta \\ 
\phi%
\end{array}%
}\right] (\tau )%
\end{array}%
\right) .
\end{equation*}%
For $\theta ,\phi \in \{\pm 1\}$, it follows from (\ref{Pkodd}) that the
diagonal Eisenstein terms vanish and that $\det \mathbf{M}=0$ for odd $n$.
Taking $n$ even and recalling that $\mathrm{Pf}(\mathbf{M})=\sqrt{\det 
\mathbf{M}}$, we recover the square of the rank one generating function (\ref%
{GnPfaff}) for $\phi =-1$ and the rank one $\sigma $-twisted generating
function (\ref{Gn_Msigma}) for $\phi =1$.

\subsection{$n$-Point Functions for $(\protect\theta ,\protect\phi )=(1,1)$.}

We consider $(\alpha ,\beta )=(0,0)$ so that $(f,g)=(1,1)$ and $(\theta
,\phi )=(1,1)$ with $\kappa =\frac{1}{2}$ (cf. Remark \ref{Remark_part_alpha}%
). We then have $wt_{h}(\psi ^{+})=$ $0$, $wt_{h}(\psi ^{-})=$ $1$ and $%
c_{h}=-2$. For $n=1$, eqn. (\ref{F_VLh_gen}) can be computed from (\ref%
{nptrec}) to give the ($x,y$ independent) result: 
\begin{eqnarray*}
G_{2,h}(1;x,y;\tau ) &=&F_{V,h}(1;(\psi ^{+},x),(\psi ^{-},y);\tau ) \\
&=&\mathrm{STr}_{V}\left( o_{h}(\psi ^{+})o_{h}(\psi
^{-})q^{L_{h}(0)+1/12}\right) +0,
\end{eqnarray*}%
where $o_{h}(v)=v(wt_{h}(v)-1)$ from (\ref{o(v)}) and recalling $%
Z_{V,h}(1;\tau )=0$. Furthermore, $o_{h}(\psi ^{+})o_{h}(\psi ^{-})=\psi
^{+}(-1)\psi ^{-}(0)$ acts as a projection operator on $V$ preserving those
Fock vectors (\ref{Fockstate}) containing an $\psi ^{+}(-1)$ operator. Hence
we find 
\begin{equation}
G_{2,h}(1;x,y;\tau )=q^{1/12}(-q^{0})\prod_{k\geq 2}(1-q^{k-1})\prod_{l\geq
1}(1-q^{l})=-\eta (\tau )^{2}.  \label{G2pt_(1,1)}
\end{equation}%
We may proceed much as before to compute the generator $G_{2n,h}$ to find:

\begin{proposition}
\label{prop_genVh_(1,1)}For $(\theta ,\phi )=(1,1)$ we have 
\begin{equation}
G_{2n,h}(1;x_{1},...,x_{n};y_{1},...,y_{n};\tau )=\det \mathbf{Q}.\eta (\tau
)^{2},  \label{G_Lh_(1,1)}
\end{equation}%
where $\mathbf{Q}$ is the $(n+1)\times (n+1)$ matrix: 
\begin{equation}
\mathbf{Q}=\left( 
\begin{array}{cccc}
P_{1}(x_{1}-y_{1},\tau ) & \ldots & P_{1}(x_{1}-y_{n},\tau ) & 1 \\ 
\vdots & \ddots &  & \vdots \\ 
P_{1}(x_{n}-y_{1},\tau ) &  & P_{1}(x_{n}-y_{n},\tau ) & 1 \\ 
1 & \ldots & 1 & 0%
\end{array}%
\right) .  \label{Qmatrix}
\end{equation}%
($P_{1}(z,\tau )$ as in (\ref{P1_dlogKdz})). Furthermore, $G_{2n,h}$ is an
analytic function in $x_{i},y_{j}$ and converges absolutely and uniformly on
compact subsets of the domain $\left\vert q\right\vert <\left\vert
q_{x_{i}-y_{j}}\right\vert <1$.
\end{proposition}

\noindent \textbf{Proof.} We prove the result by induction. For $n=1$ we
obtain the result from (\ref{G2pt_(1,1)})$.$ Assuming the result for $n-1$,
we apply the Zhu recursive formula (\ref{nptrec}) to find%
\begin{eqnarray*}
&&G_{2n,h}(1;x_{1},...,x_{n};y_{1},...,y_{n};\tau ) \\
&=&\mathrm{STr}_{V}\left( o(\psi ^{+})Y(q_{y_{1}}^{L(0)}\psi
^{-},q_{y_{1}})\ldots Y(q_{x_{n}}^{L(0)}\psi
^{+},q_{x_{n}})Y(q_{y_{n}}^{L(0)}\psi ^{-},q_{y_{n}})q^{L(0)+1/12}\right) \\
&&+\sum\limits_{r=1}^{n}(-1)^{r-1}\mathbf{Q}(1,r)\mathrm{\det }\mathbf{\hat{Q%
}}.\eta (\tau )^{2},
\end{eqnarray*}%
where $\mathbf{\hat{Q}}$ denotes the matrix found from $\mathbf{Q}$ by
deleting row $1$ and column $r$. Next note from Lemma \ref{Lemma_Perm_Period}
(ii) and (\ref{psiplus_minus_comm}) that $G_{2n,h}$ vanishes for $%
x_{1}=x_{2} $ so that 
\begin{eqnarray*}
&&\mathrm{STr}_{V}\left( o(\psi ^{+})Y(q_{y_{1}}^{L(0)}\psi
^{-},q_{y_{1}})\ldots Y(q_{y_{n}}^{L(0)}\psi
^{-},q_{y_{n}})q^{L(0)+1/12}\right) \\
&=&-\sum\limits_{r=2}^{n}(-1)^{r}\mathbf{Q}(2,r)\mathrm{\det }\mathbf{\hat{Q}%
}.\eta (\tau )^{2}.
\end{eqnarray*}%
Hence we find $G_{2n,h}$ is given by

\begin{equation*}
\sum\limits_{r=1}^{n}(-1)^{r-1}(\mathbf{Q}(1,r)-\mathbf{Q}(2,r))\mathrm{\det 
}(\mathbf{\hat{Q}})\eta (\tau )^{2}=\det \mathbf{Q.}\eta (\tau )^{2},
\end{equation*}%
on evaluating $\det \mathbf{Q}$ after subtracting row 2 from row 1. $\square 
$

\medskip We may similarly obtain a determinant formula for all $n$-point
functions along the same lines as Propositions \ref{Prop_Ferm_1_npt} and \ref%
{prop_npt_Vh_g}.

\bigskip

\subsection{\label{Boson}Bosonization}

As is well known, the rank two fermion VOSA $V$ can be constructed as a rank
one bosonic $\mathbb{Z}$-lattice VOSA. $V$ is decomposed in terms of the
Heisenberg subVOA $M$ generated by the boson $a$ of (\ref{adef}) and its
irreducible modules $M\otimes e^{m}$ for $a(0)$ eigenvalue $m\in \mathbb{Z}$
(cf. \cite{Ka}). In particular, the partition function $Z_{V,h}(f;\tau )$
and the generating function $G_{n,h}$ can be computed in this bosonic
decomposition using the results of ref. \cite{MT1}, leading to the Jacobi
triple product formula and Fay's trisecant identity (for elliptic functions)
respectively. We also describe a further new generalization of Fay's
trisecant identity for elliptic functions.

The highest weight lattice vector for the irreducible module $M\otimes e^{m}$
is 
\begin{equation*}
\mathbf{1}\otimes e^{m}=\left\{ 
\begin{array}{c}
\psi ^{+}(-m)\psi ^{+}(1-m)...\psi ^{+}(-1).\mathbf{1},\quad m>0, \\ 
\psi ^{-}(m)\psi ^{-}(1+m)...\psi ^{-}(-1).\mathbf{1},\quad m<0.%
\end{array}%
\right.
\end{equation*}%
Then the partition function is 
\begin{eqnarray}
Z_{V,h}(f;\tau ) &=&Z_{V}\left[ 
\begin{array}{c}
f \\ 
g%
\end{array}%
\right] (\tau )=\sum_{m\in \mathbb{Z}}(-1)^{m}e^{2\pi im\alpha }\mathrm{Tr}%
_{M\otimes e^{m}}(q^{L(0)+\kappa ^{2}/2-\kappa m-1/24})  \notag \\
&=&\frac{e^{2\pi i(\alpha +1/2)(\beta +1/2)}}{\eta (\tau )}\vartheta \left[ 
\begin{array}{c}
-\beta +\frac{1}{2} \\ 
\alpha +\frac{1}{2}%
\end{array}%
\right] (0,\tau ),  \label{Z_lattice}
\end{eqnarray}%
in terms of the theta series (\ref{thetaab}). Comparing to (\ref{Zparth}) we
obtain the standard Jacobi triple product formula.

\medskip We can also compute the generating function $G_{n,h}$ (and hence
all $n$-point functions) in the bosonic setting based on results of ref. 
\cite{MT1}. We illustrate this with the 2-point function generator (\ref%
{F_VLh_gen}). Recall from (\ref{FMonepointgen11}) and (\ref{G2pt_(1,1)})
that 
\begin{equation}
G_{2,h}(f;x;y;\tau )=\left\{ 
\begin{array}{cc}
P_{1}\left[ 
\begin{array}{c}
\theta \\ 
\phi%
\end{array}%
\right] (x-y,\tau )Z_{V,h}(f;\tau ), & (\theta ,\phi )\neq (1,1), \\ 
-\eta (\tau )^{2}, & (\theta ,\phi )=(1,1).%
\end{array}%
\right.  \label{G2h}
\end{equation}%
In the bosonic language we obtain: 
\begin{eqnarray*}
G_{2,h}(f;x;y;\tau ) &=&\sum_{m\in \mathbb{Z}}(-1)^{m}e^{2\pi im\alpha
}F_{M\otimes e^{m},h}(1;(\mathbf{1}\otimes e^{+1},x),(\mathbf{1}\otimes
e^{-1},y);\tau ) \\
&=&\sum_{m\in \mathbb{Z}}(-1)^{m}e^{2\pi im\alpha }\exp (-\kappa
(x-y))q^{\kappa ^{2}/2-\kappa m}. \\
&&F_{M\otimes e^{m}}(1;(\mathbf{1}\otimes e^{+1},x),(\mathbf{1}\otimes
e^{-1},y);\tau ),
\end{eqnarray*}%
noting that $Y(q_{z}^{L_{h}(0)}e^{\pm 1},q_{z})=\exp (\mp \kappa
z)Y(q_{z}^{L(0)}e^{\pm 1},q_{z})$. Using Propositions 4 and 5 of ref. \cite%
{MT1} we obtain%
\begin{equation*}
F_{M\otimes e^{m}}(1;(\mathbf{1}\otimes e^{+1},x),(\mathbf{1}\otimes
e^{-1},y);\tau )=\frac{q^{m^{2}/2}}{\eta (\tau )}\frac{\exp (m(x-y))}{%
K(x-y,\tau )},
\end{equation*}%
where $K$ is the prime form (\ref{Primeform}). Altogether, it follows that%
\begin{equation*}
G_{2,h}(f;x;y;\tau )=\frac{e^{2\pi i(\alpha +1/2)(\beta +1/2)}}{\eta (\tau )}%
\frac{\vartheta \left[ 
\begin{array}{c}
-\beta +\frac{1}{2} \\ 
\alpha +\frac{1}{2}%
\end{array}%
\right] (x-y,\tau )}{K(x-y,\tau )}.
\end{equation*}%
Comparing with (\ref{G2h}) we confirm the identities (\ref{P1uv_theta}) for $%
(\theta ,\phi )\neq (1,1)$ and (\ref{Ktheta}) for $(\alpha ,\beta )=(0,0)$,
i.e. $(\theta ,\phi )=(1,1)$.

\medskip In a similar fashion we can compute the general generating function 
$G_{2n,h} $ in the bosonic setting to obtain:

\begin{proposition}
\label{Proposition_G2n_boson}%
\begin{gather*}
G_{2n,h}(f;x_{1},...,x_{n};y_{1},...,y_{n};\tau )=\frac{e^{2\pi i(\alpha
+1/2)(\beta +1/2)}}{\eta (\tau )}\vartheta \left[ 
\begin{array}{c}
-\beta +\frac{1}{2} \\ 
\alpha +\frac{1}{2}%
\end{array}%
\right] (\sum_{i=1}^{n}(x_{i}-y_{i}),\tau ). \\
\frac{\prod\limits_{1\leq i<j\leq n}K(x_{i}-x_{j},\tau )K(y_{i}-y_{j},\tau )%
}{\prod\limits_{1\leq i,j\leq n}K(x_{i}-y_{j},\tau )}.\quad \square
\end{gather*}
\end{proposition}

Comparing this to Proposition \ref{prop_genVh} for $(\theta ,\phi )\neq
(1,1) $ and Proposition \ref{prop_genVh_(1,1)} for $(\theta ,\phi )=(1,1)$
we obtain the elliptic function version of Fay's Generalized Trisecant
Identity \cite{Fa}:

\begin{corollary}
For $(\theta ,\phi )\neq (1,1)$ we have 
\begin{equation}
\mathrm{\det }(\mathbf{P})=\frac{\vartheta \left[ 
\begin{array}{c}
-\beta +\frac{1}{2} \\ 
\alpha +\frac{1}{2}%
\end{array}%
\right] (\sum_{i=1}^{n}(x_{i}-y_{i}),\tau )}{\vartheta \left[ 
\begin{array}{c}
-\beta +\frac{1}{2} \\ 
\alpha +\frac{1}{2}%
\end{array}%
\right] (0,\tau )}\frac{\prod\limits_{1\leq i<j\leq n}K(x_{i}-x_{j},\tau
)K(y_{i}-y_{j},\tau )}{\prod\limits_{1\leq i,j\leq n}K(x_{i}-y_{j},\tau )},
\label{Fay_trisecant}
\end{equation}%
with $\mathbf{P}$ as in (\ref{Pmatrix}). For $(\theta ,\phi )=(1,1)$, 
\begin{equation}
\mathrm{\det }(\mathbf{Q})=-\frac{K(\sum_{i=1}^{n}(x_{i}-y_{i}),\tau
)\prod\limits_{1\leq i<j\leq n}K(x_{i}-x_{j},\tau )K(y_{i}-y_{j},\tau )}{%
\prod\limits_{1\leq i,j\leq n}K(x_{i}-y_{j},\tau )},  \label{K_secant}
\end{equation}%
where $\mathbf{Q}$ is as in (\ref{Qmatrix}). $\square $
\end{corollary}

We may generalize these identities using Propositions 4 and 5 of \cite{MT1}
again to consider the general lattice $n$-point function:

\begin{proposition}
\label{Proposition_Flattice_boson}For integers $m_{i},n_{j}\geq 0$
satisfying 
\begin{equation*}
\sum_{i=1}^{r}m_{i}=\sum_{j=1}^{s}n_{j},
\end{equation*}
we have 
\begin{eqnarray*}
&&F_{V}(f;(\mathbf{1\otimes }e^{m_{1}},x_{1}),...(\mathbf{1\otimes }%
e^{m_{r}},x_{r}),(\mathbf{1\otimes }e^{-n_{1}},y_{1}),...(\mathbf{1\otimes }%
e^{-n_{s}},y_{s});\tau ) \\
&=&\frac{e^{2\pi i(\alpha +1/2)(\beta +1/2)}}{\eta (\tau )}\vartheta \left[ 
\begin{array}{c}
-\beta +\frac{1}{2} \\ 
\alpha +\frac{1}{2}%
\end{array}%
\right] (\sum_{i=1}^{r}m_{i}x_{i}-\sum_{j=1}^{s}n_{j}y_{j},\tau ). \\
&&\frac{\prod\limits_{1\leq i<k\leq r}K(x_{i}-x_{k},\tau
)^{m_{i}m_{k}}\prod\limits_{1\leq j<l\leq s}K(y_{j}-y_{l},\tau )^{n_{j}n_{l}}%
}{\prod\limits_{1\leq i\leq r,1\leq j\leq s}K(x_{i}-y_{j},\tau )^{m_{i}n_{j}}%
}.\quad \square
\end{eqnarray*}
\end{proposition}

Comparing this to Proposition \ref{prop_npt_Vh_g} we obtain a new elliptic
generalization of Fay's Trisecant Identity:

\begin{corollary}
For $(\theta ,\phi )\neq (1,1)$ we have%
\begin{eqnarray*}
\mathrm{\det }(\mathbf{M}) &=&\frac{\vartheta \left[ 
\begin{array}{c}
-\beta +\frac{1}{2} \\ 
\alpha +\frac{1}{2}%
\end{array}%
\right] (\sum_{i=1}^{r}m_{i}x_{i}-\sum_{j=1}^{s}n_{j}y_{j},\tau )}{\vartheta %
\left[ 
\begin{array}{c}
-\beta +\frac{1}{2} \\ 
\alpha +\frac{1}{2}%
\end{array}%
\right] (0,\tau )}. \\
&&\frac{\prod\limits_{1\leq i<k\leq r}K(x_{i}-x_{k},\tau
)^{m_{i}m_{k}}\prod\limits_{1\leq j<l\leq s}K(y_{j}-y_{l},\tau )^{n_{j}n_{l}}%
}{\prod\limits_{1\leq i\leq r,1\leq j\leq s}K(x_{i}-y_{j},\tau )^{m_{i}n_{j}}%
},
\end{eqnarray*}%
where $\mathbf{M}$ is the block matrix%
\begin{equation*}
\mathbf{M}=\left( 
\begin{array}{ccc}
\mathbf{D}^{(11)} & \ldots & \mathbf{D}^{(1s)} \\ 
\vdots & \ddots & \vdots \\ 
\mathbf{D}^{(r1)} & \ldots & \mathbf{D}^{(rs)}%
\end{array}%
\right) ,
\end{equation*}%
with $\mathbf{D}^{(ab)}$ the $m_{a}\times n_{b}$ matrix 
\begin{equation*}
\mathbf{D}^{(ab)}(i,j)=D\left[ 
\begin{array}{c}
\theta \\ 
\phi%
\end{array}%
\right] (i,j,\tau ,x_{a}-y_{b}),\quad (1\leq i\leq m_{a},1\leq j\leq n_{b}),
\end{equation*}%
for $1\leq a\leq r$ and $1\leq b\leq s$. $\square $
\end{corollary}

A similar identity for $(\theta ,\phi )=(1,1)$ generalizing (\ref{K_secant})
can also be described.

\subsection{Modular Properties of $n$-Point Functions}

In this section we consider the modular properties of all $n$-point
functions for the rank two fermion VOSA. Despite the fact the twisted
sectors are neither rational or $C_{2}$-cofinite we obtain modular
properties similar to those found in \cite{Z}, \cite{DZ1}, \cite{DZ2}. It is
convenient to employ the twisted $n$-point function formalism to describe
these modular properties. We firstly consider the partition function $Z_{V}%
\left[ 
\begin{array}{c}
f \\ 
g%
\end{array}%
\right] (\tau )$ and define a group action for $\gamma =\left( 
\begin{array}{cc}
a & b \\ 
c & d%
\end{array}%
\right) \in SL(2,\mathbb{Z})$ as follows: 
\begin{equation}
\left. Z_{V}\left[ 
\begin{array}{c}
f \\ 
g%
\end{array}%
\right] \right\vert \gamma (\tau) =Z_{V}\left( \gamma .\left[ 
\begin{array}{c}
f \\ 
g%
\end{array}%
\right] \right) (\gamma .\tau ),  \label{ZV_gamma}
\end{equation}%
with 
\begin{equation}
\gamma .\left[ 
\begin{array}{c}
f \\ 
g%
\end{array}%
\right] =\left[ 
\begin{array}{c}
f^{a}g^{b} \\ 
f^{c}g^{d}%
\end{array}%
\right] ,  \label{gamma_fg}
\end{equation}%
and $\gamma .\tau $ as in (\ref{mod_z_tau}).

\begin{remark}
\label{Remark_gamma_left}(i) (\ref{gamma_fg}) is equivalent to left matrix
multiplication on $\alpha ,\beta $ 
\begin{equation*}
\gamma \left( 
\begin{array}{c}
\alpha \\ 
\beta%
\end{array}%
\right) =\left( 
\begin{array}{c}
a\alpha +b\beta \\ 
c\alpha +d\beta%
\end{array}%
\right) .
\end{equation*}

(ii) In terms of the shifted VOA formalism, (\ref{ZV_gamma}) reads 
\begin{equation*}
\left. Z_{V,h}(f;\tau )\right\vert \gamma = Z_{V,\gamma
.h}(f^{a}g^{b};\gamma .\tau ),
\end{equation*}%
with $\gamma .h=(\gamma .\beta +\frac{1}{2})a=((c\alpha +d\beta )+\frac{1}{2}%
)a,$ recalling (\ref{h_kappa}) and (\ref{beta_kappa}).
\end{remark}

For $SL(2,\mathbb{Z})$ generators $S=\left( 
\begin{array}{cc}
0 & 1 \\ 
-1 & 0%
\end{array}%
\right) $ and $T=\left( 
\begin{array}{cc}
1 & 1 \\ 
0 & 1%
\end{array}%
\right) $ we can use the theta function modular transformation properties (%
\ref{S_Theta}) and (\ref{T_Theta}) and thereby find from (\ref{Z_lattice})
that 
\begin{eqnarray}
\left. Z_{V}\left[ 
\begin{array}{c}
f \\ 
g%
\end{array}%
\right] \right\vert S(\tau ) &=&\varepsilon _{S}\left[ 
\begin{array}{c}
f \\ 
g%
\end{array}%
\right] Z_{V}\left[ 
\begin{array}{c}
f \\ 
g%
\end{array}%
\right] (\tau ),  \label{ZV_S} \\
\left. Z_{V}\left[ 
\begin{array}{c}
f \\ 
g%
\end{array}%
\right] \right\vert T(\tau ) &=&\varepsilon _{T}\left[ 
\begin{array}{c}
f \\ 
g%
\end{array}%
\right] Z_{V}\left[ 
\begin{array}{c}
f \\ 
g%
\end{array}%
\right] (\tau ),  \label{ZV_T}
\end{eqnarray}%
where 
\begin{eqnarray}
\varepsilon _{S}\left[ 
\begin{array}{c}
f \\ 
g%
\end{array}%
\right]  &=&\exp (2\pi i(\frac{1}{2}+\beta )(\frac{1}{2}-\alpha )), \\
\varepsilon _{T}\left[ 
\begin{array}{c}
f \\ 
g%
\end{array}%
\right]  &=&\exp (\pi i(\beta (\beta +1)+\frac{1}{6})).  \label{eps_ST}
\end{eqnarray}%
One can check that the relations $(ST)^{3}=-S^{2}=1$ are satisfied so that $%
Z_{V}\left[ 
\begin{array}{c}
f \\ 
g%
\end{array}%
\right] (\tau )$ is modular invariant as follows:

\begin{proposition}
\label{Proposition_ZV_modinv}The partition function transforms under $\gamma
\in SL(2, \mathbb{Z})$  with multiplier $\varepsilon _{\gamma }\left[ 
\begin{array}{c}
f \\ 
g%
\end{array}%
\right] \in U(1)$ where 
\begin{equation*}
\left. Z_{V}\left[ 
\begin{array}{c}
f \\ 
g%
\end{array}%
\right] \right\vert \gamma (\tau)=\varepsilon _{\gamma }\left[ 
\begin{array}{c}
f \\ 
g%
\end{array}%
\right] Z_{V}\left[ 
\begin{array}{c}
f \\ 
g%
\end{array}%
\right] (\tau ),
\end{equation*}
with $\varepsilon _{\gamma }\left[ 
\begin{array}{c}
f \\ 
g%
\end{array}%
\right] $ generated from $\varepsilon _{S}\left[ 
\begin{array}{c}
f \\ 
g%
\end{array}%
\right] $ and $\varepsilon _{T}\left[ 
\begin{array}{c}
f \\ 
g%
\end{array}%
\right] $. $\square $
\end{proposition}

In order to discuss the modular properties of $n$-point functions we first
define the left $SL(2,\mathbb{Z})$ action%
\begin{equation}
\left. F_{V}\left[ 
\begin{array}{c}
f \\ 
g%
\end{array}%
\right] ((v_{1},z_{1})...,(v_{n},z_{n});\tau )\right\vert \gamma
=F_{V}\left( \gamma .\left[ 
\begin{array}{c}
f \\ 
g%
\end{array}%
\right] \right) ((v_{1},\gamma .z_{1})...,(v_{n},\gamma .z_{n});\gamma .\tau
),  \label{FV_gamma}
\end{equation}%
and $\gamma .z$ as in (\ref{mod_z_tau}). It is sufficient to consider the
generating function:

\begin{proposition}
\label{Proposition_Gen_mod}The generating function $G_{2n}\left[ 
\begin{array}{c}
f \\ 
g%
\end{array}%
\right] $ transforms under $\gamma \in SL(2, \mathbb{Z})$ with weight $n$
and multiplier $\varepsilon _{\gamma }\left[ 
\begin{array}{c}
f \\ 
g%
\end{array}%
\right] $, that is 
\begin{equation*}
\left. G_{2n}\left[ 
\begin{array}{c}
f \\ 
g%
\end{array}%
\right] (x_{1}...x_{n};y_{1}...y_{n};\tau )\right\vert \gamma =(c\tau
+d)^{n}\varepsilon _{\gamma }\left[ 
\begin{array}{c}
f \\ 
g%
\end{array}%
\right] G_{2n}\left[ 
\begin{array}{c}
f \\ 
g%
\end{array}%
\right] (x_{1}...x_{n};y_{1}...y_{n};\tau ).
\end{equation*}
\end{proposition}

\noindent \textbf{Proof. }For $(\theta ,\phi )\neq (1,1)$ we have 
\begin{equation*}
G_{2n}\left[ 
\begin{array}{c}
f \\ 
g%
\end{array}%
\right] =\mathrm{\det }(\mathbf{P})\;Z_{V}\left[ 
\begin{array}{c}
f \\ 
g%
\end{array}%
\right] (\tau ),
\end{equation*}%
from Proposition \ref{prop_genVh}. From Proposition \ref%
{Proposition_Pk_Modular} we have $P_{1}\left[ 
\begin{array}{c}
\theta \\ 
\phi%
\end{array}%
\right] (\gamma .z,\gamma .\tau )=(c\tau +d)P_{1}\left[ 
\begin{array}{c}
\theta \\ 
\phi%
\end{array}%
\right] (z,\tau )$. Hence using Proposition \ref{Proposition_ZV_modinv} the
result follows.

For $(\theta ,\phi )=(1,1)$ we have%
\begin{equation*}
G_{2n}\left[ 
\begin{array}{c}
f \\ 
g%
\end{array}%
\right] =\det \mathbf{Q}.\eta (\tau )^{2},
\end{equation*}%
from Proposition \ref{prop_genVh_(1,1)} with $\mathbf{Q}$ as in (\ref%
{Qmatrix}). From (\ref{gammaEn}) and (\ref{gammaE2}) it follows that $%
P_{1}(z,\tau )$ is quasi-modular: 
\begin{equation*}
P_{1}(\gamma .z,\gamma .\tau )=(c\tau +d)P_{1}(z,\tau )+\frac{c}{2\pi i}z.
\end{equation*}%
However, $\det \mathbf{Q}$ is modular of weight $n-1$ as follows. Subtract
row 1 from rows $2\ldots n$ and then subtract col 1 from cols $2\ldots n$ to
find $\det \mathbf{Q=}\det \mathbf{R}$ where for $2\leq i,j\leq n$ 
\begin{equation*}
R(i,j)=P_{1}(x_{i}-y_{j},\tau )+P_{1}(x_{1}-y_{1},\tau
)-P_{1}(x_{i}-y_{1},\tau )-P_{1}(x_{1}-y_{j},\tau ),
\end{equation*}%
which is modular of weight $1$. Hence the result follows. $\square $

\medskip The modular transformation properties for an arbitrary $n$-point
function follows by appropriately expanding the generating function as
before to find $n$-point functions for the Fock basis described in
Proposition \ref{prop_npt_Vh_g}. We thus find

\begin{proposition}
\label{Proposition_npt_mod}For $n$ vectors $v_{a}$ of $wt[v_{a}]$, $a=1,
\ldots, n$, the $n$-point function transforms under $\gamma \in SL(2, 
\mathbb{Z})$  with weight $K=\sum_{a}wt[v_{a}]$ and multiplier $\varepsilon
_{\gamma }\left[ 
\begin{array}{c}
f \\ 
g%
\end{array}%
\right] $: 
\begin{multline*}
\left. F_{V}\left[ 
\begin{array}{c}
f \\ 
g%
\end{array}%
\right] ((v_{1},z_{1}),\ldots ,(v_{n},z_{n});\tau )\right\vert \gamma = \\
(c\tau +d)^{K}\varepsilon _{\gamma }\left[ 
\begin{array}{c}
f \\ 
g%
\end{array}%
\right] F_{V}\left[ 
\begin{array}{c}
f \\ 
g%
\end{array}%
\right] ((v_{1},z_{1}),\ldots ,(v_{n},z_{n});\tau).
\end{multline*}
\end{proposition}

This result is a natural generalization for continuous orbifolds of the rank
two fermion VOSA of Zhu's Theorem 5.3.2 for $C_{2}$-cofinite VOAs \cite{Z}.

\section{Appendix A: Parity and Supertraces}


A vertex operator $Y(a,z)$ has parity $p(a)\in \{0,1\}$ if all its modes $%
a(n)$ have parity $p(a)$. Two operators $A,B$ on $V$ of parity $p(A),p(B)$
have commutator defined by 
\begin{eqnarray*}
\lbrack A,B] &=&AB-p(A,B)BA, \\
p(A,B) &=&(-1)^{p(A)p(B)}.
\end{eqnarray*}%
The commutator clearly obeys: 
\begin{equation*}
\lbrack A,B]=-p(A,B)[B,A],
\end{equation*}%
and for $B_{1}\ldots B_{n}$ of parity $p(B_{1}),\ldots p(B_{n})$
respectively we have 
\begin{eqnarray}
&&[A,B_{1}\ldots B_{n}]  \notag \\
&=&\sum\limits_{r=1}^{n}p(A,B_{1}\ldots B_{r-1})B_{1}\ldots
B_{r-1}[A,B_{r}]B_{r+1}\ldots B_{n},  \label{AB1toBncom}
\end{eqnarray}%
where%
\begin{equation}
p(A,B_{1}\ldots B_{r-1})=\left\{ 
\begin{array}{cc}
1\text{ } & \text{for }r=1 \\ 
(-1)^{p(A)[p(B_{1})+...+p(B_{r-1})]} & \text{ for }r>1%
\end{array}%
\right. .  \label{parityAB}
\end{equation}

Let\ $V_{\alpha }=\bigoplus\limits_{r\geq r_{0}}V_{\alpha ,r}$ denoted the
decomposition of $V_{\alpha }$ into $L(0)$ homogeneous spaces where $r_{0}$
is the lowest $L(0)$ degree. We assume that $\dim V_{\alpha ,r}$ is finite
for each $r,\alpha $. We define the Supertrace of an operator $A$ by: 
\begin{eqnarray*}
\mathrm{STr}(Aq^{L(0)}) &=&Tr(\sigma Aq^{L(0)}) \\
&=&Tr_{V_{\bar{0}}}(Aq^{L(0)})-Tr_{V_{\bar{1}}}(Aq^{L(0)}) \\
&=&\sum\limits_{r\geq r_{0}}q^{r}[Tr_{V_{\bar{0},r}}(A)-Tr_{V_{\bar{1}%
,r}}(A)].
\end{eqnarray*}%
Clearly the supertrace is zero if $A$ has odd parity. We then note the
following:

\begin{lemma}
\label{lemmacommk} Suppose that $A$ is an operator on $V$ of parity $p(A)$
such that $A:V_{\alpha ,r}\rightarrow V_{\alpha +p(A),r+s}$ for some real $s$%
. Then for any operator $B$ we have: 
\begin{equation*}
\mathrm{STr}(ABq^{L(0)})=q^{s}p(A,B)\;\mathrm{STr}(BAq^{L(0)}).\quad \square
\end{equation*}
\end{lemma}

Using (\ref{bnVr}) we find

\begin{corollary}
\label{corrv(k)}For $v$ homogeneous of weight $wt(v)$ then 
\begin{equation}
\mathrm{STr}(v(k)\;B\;q^{L(0)})=p(v,B)q^{wt(v)-k-1}\;\mathrm{STr}%
(Bv(k)q^{L(0)}).\quad \square  \label{corrv(k)formula}
\end{equation}
\end{corollary}

We also have

\begin{corollary}
\label{corrcomm} If $A:V_{\alpha ,r}\rightarrow V_{\alpha +p(A),r}$ then for
any operator $B$ we have 
\begin{equation}
\mathrm{STr}([A,B]q^{L(0)})=0.\quad \square  \label{StrAB}
\end{equation}
\end{corollary}

\end{document}